\documentclass[11pt]{amsart}

\usepackage{amsmath}
\usepackage{fullpage}
\usepackage{xspace}
\usepackage[psamsfonts]{amssymb}
\usepackage{graphicx,color}
\usepackage{hyperref}
\usepackage{graphicx}
\usepackage{amsmath}%
\usepackage{amsthm}%
\usepackage{amscd}
\usepackage{amsfonts}%
\usepackage{amssymb}%
\usepackage{graphicx}
\usepackage{stmaryrd}
\usepackage{mathrsfs}
\usepackage{tikz}
\usetikzlibrary{matrix,arrows}
\usepackage{tikz-cd}
\usepackage[all]{xy}
\usepackage{bbm} % para mathbb en números
\usepackage{longtable}
\usepackage{comment}

\newcounter{Main}

\theoremstyle{plain} % Pour tous ce qui ressemble à des théorèmes.

\newtheorem{MainThm}[Main]{Theorem} % Théorèmes principaux
\newtheorem{MainConj}[Main]{Conjecture} % Corollaires principaux.
\newtheorem{theorem}{Theorem}[section]

\newtheorem{corollary}[theorem]{Corollary}

\newtheorem{definition}{Definition}
\newtheorem{example}{Example}

\newtheorem{lemma}[theorem]{Lemma}

\newtheorem{proposition}[theorem]{Proposition}
\theoremstyle{remark}
\newtheorem{remark}[theorem]{Remark}

\numberwithin{equation}{section}

\newcommand{\X}{X}
\newcommand{\Xinf}{X_{\infty}}
\newcommand{\OS}{\mathcal{O}_\mathcal{S}}
\newcommand{\rk}{\mathrm{rk}}
\newcommand{\St}{\mathrm{St}}
\newcommand{\ZG}{\mathbb{Z}[G]}

\newcommand\ssbmatrix[4]{\textnormal{\scriptsize$\left(\begin{array}{cc}#1&#2\\#3&#4\end{array}\right)$\normalsize}}

\usepackage{xcolor}
\usepackage[all]{xy}
\newcommand\rquot[2]{
  \mathchoice
  {% \displaystyle
    \text{\raise0.5ex\hbox{$#1$}\big/\lower0.5ex\hbox{$#2$}}%
  }
  {% \textstyle
    #1\,/\,#2
  }
  {% \scriptstyle
    #1\,/\,#2
  }
  {% \scriptscriptstyle
    #1\,/\,#2
  }
}

\newcommand\lrquot[3]{
  \mathchoice
  {% \displaystyle
    \text{\lower0.5ex\hbox{$#1$}\big\backslash\raise0.5ex\hbox{$#2$\!}\big/
      \lower0.5ex\hbox{\!\!$#3$}}%
  }
  {% \textstyle
    #1\,\backslash\,#2\,/\,#3
  }
  {% \scriptstyle
    #1\,\backslash\,#2\,/\,#3
  }
  {% \scriptscriptstyle
    #1\,\backslash\,#2\,/\,#3
  }
}

\newcommand\lquot[2]{
  \mathchoice
  {% \displaystyle
    \text{\lower0.5ex\hbox{$#1$}\big\backslash\raise0.5ex\hbox{$#2$}}%
  }
  {% \textstyle
    #1\,\backslash\,#2
  }
  {% \scriptstyle
    #1\,\backslash\,#2
  }
  {% \scriptscriptstyle
    #1\,\backslash\,#2
  }
}

\usepackage{fontawesome}
  \DeclareFontFamily{U}{wncy}{}
    \DeclareFontShape{U}{wncy}{m}{n}{<->wncyr10}{}
    \DeclareSymbolFont{mcy}{U}{wncy}{m}{n}
    \DeclareMathSymbol{\Sha}{\mathord}{mcy}{"58}

\begin{document}
\title[Relative homology of $\mathcal{S}$-arithmetic groups of $\mathrm{SL}_2$]{Relative homology of $\mathcal{S}$-arithmetic subgroups of $\mathrm{SL}_2$ in positive characteristic}
\author{Matias Alvarado}
\address{Instituto de Matem\'aticas, Universidad de Talca, Talca, Chile.}
\email[M. Alvarado]{matias.alvarado@utalca.cl}
\author{Claudio Bravo}
\address{Instituto de Matem\'aticas, Universidad de Talca, Talca, Chile.}
\email[C. Bravo]{claudio.bravo@utalca.cl}
%\email[P. P\'erez-pi\~na]{paperez15@uc.cl }%
%\thanks{Supported by ANID Doctorado Nacional 21200910}
\date{\today}

\begin{abstract}
%Let $k$ be a global field of positive characteristic. We describe the relative homology groups of certain $S$-arithmetic subgroups $G$ of $\mathrm{SL}_2(k)$, modulo a set of representatives for the conjugacy classes of their maximal unipotent subgroups, in terms of the Euler-Poincar\'e characteristic. We also introduce the notion of relative homological dimension and prove that its value for $G$ equals the cardinality of $S$.
Let $k$ be a global field of positive characteristic $p$, let $\mathcal{S}$ be a non-empty finite set of places of $k$, and let $G$ be an $\mathcal{S}$-arithmetic subgroup of $\mathrm{SL}_2(k)$ whose torsion is $p$-primary. We study the relative homology of $G$ modulo two natural families of subgroups, which are the stabilizers of the cusps, arising from the action of $G$ on the product of the Bruhat-Tits trees at the places of $\mathcal{S}$, and a set of representatives of the conjugacy classes of maximal unipotent subgroups of $G$. 
We prove that the Steinberg module of $G$ defined by the cusp stabilizers has projective and flat dimension $\#\mathcal{S}-1$.
In particular, and in contrast with
the classical case $\#\mathcal{S}=1$, it is not projective as soon as $\#\mathcal{S}\geq 2$.
Consequently, the relative homological dimension of $G$ modulo the cusp stabilizers equals $\#\mathcal{S}$, and we prove that the same holds modulo the maximal unipotent subgroups whenever $p$ is invertible in the coefficient ring. 
We express
the alternating sum of the ranks of the corresponding relative homology groups in terms of the Euler--Poincar\'e characteristic of $G$. For maximal unipotent subgroups and $\#\mathcal{S}\geq 2$, a correction term appears which, for trivial coefficients, equals
the number of cusps of $G$. %In particular, a conjecture of Serre on the relative homology of $\mathcal{S}$-arithmetic groups modulo their maximal unipotent subgroups does not hold as stated.
Finally, we show that the homological dimension of $G$ over any ring in which $p$ is invertible equals $\#\mathcal{S}$.\\
\textbf{MSC codes:} primary 20J05, 20G30, 11F75; secondary 20E42, 11R58.\\
\textbf{Keywords:} arithmetic groups, function fields, Bruhat-Tits buildings, relative homology, maximal unipotent subgroups, homological dimension.
\end{abstract}
\maketitle
\setcounter{tocdepth}{1}

%\tableofcontents

\section{Introduction}\label{introduction}

\begin{comment}
Let $\mathtt{G}$ be a reductive group defined over a field $\mathtt{K}$ endowed with a discrete valuation. One of the most efficient way to study $\mathtt{G}(\mathtt{K})$, and its discrete subgroup, is through their actions on the Bruhat-Tits building $X(\mathtt{G},\mathtt{K})$, which is a CW-complex introduced by Bruhat and Tits in \cite{BT1,BT2}.
This is a contractible simplicial complex, on which $\mathtt{G}(\mathtt{K})$ acts by simplicial automorphisms.
The building is also a natural tool to compute the (co)homology of discrete subgroups of $\mathtt{G}(\mathtt{K})$.
In the function field context, this point of view has been fruitful exploited by Serre \cite{SerreTrees}, Soul\'e \cite{soule}, Margaux \cite{margaux}, and Harder \cite{harder}, among many others. 
%This is the point of view we adopt in the present work.
\end{comment}

Let $\mathtt{G}$ be a reductive group defined over a field $\mathtt{K}$ endowed with a discrete valuation. One of the most efficient ways to study $\mathtt{G}(\mathtt{K})$ and its discrete subgroups is through their action on the Bruhat-Tits building $\X=\X(\mathtt{G},\mathtt{K})$, introduced by Bruhat and Tits in \cite{BT1,BT2}.
This is a contractible complex on which $\mathtt{G}(\mathtt{K})$ acts
by simplicial automorphisms, which makes it particularly well suited to the study of the (co)homology of discrete subgroups $\Gamma\subseteq \mathtt{G}(\mathtt{K})$.
Indeed, since $\X$ is contractible, the equivariant spectral sequence attached to the action of $\Gamma$ computes the homology $H_*(\Gamma,M)$ in terms of the homology of the cell stabilizers and of the combinatorics of the quotient $\Gamma\backslash \X$ (cf. \cite[Ch.~VII]{Brown-Cohomology}). 
%The same action also provides a natural definition of the Euler-Poincar\'e characteristic of $\Gamma$ (cf. \cite[??]{SerreCohomology}).

In the function field setting, this point of view has been fruitfully
exploited by many authors.
Serre, for instance,  described the quotient of the Bruhat-Tits tree by arithmetic subgroups of $\mathrm{SL}_2$ and deduced from it the structure of their homology modulo the cusps \cite[Ch.~II]{SerreTrees}.

Soul\'e \cite{soule} showed that, for Chevalley group schemes over polynomial rings, the quotient is a sector chamber. Equivalently, this holds for split semisimple simply connected $\mathtt{K}$-groups, where a sector chamber is a fundamental domain for the action of a Weyl group on a tessellated finite-dimensional vector space.

This result has been extended to isotrivial groups by Margaux in \cite{margaux}.
All these results yield amalgam decompositions and explicit homology computations for the groups involved.
Harder used the aforementioned approach to study the cohomology of $S$-arithmetic groups over function fields and its relation
with automorphic forms \cite{harder}.

Finiteness properties of $S$-arithmetic groups in positive characteristic are considerably more subtle than in the number field case (cf. \cite{Gandini}).
%Stuhler showed in \cite{Stuhler} that $\mathrm{SL}_2(\OS)$ is of type $F_{\sharp S-1}$ but not of type $F_{\sharp S}$, and that it has cohomology groups of infinite dimension.
Stuhler studied in \cite{Stuhler} the finiteness properties of $\mathrm{SL}_2(\OS)$, namely the properties $FP_n$, which generalize the notions of being finitely generated and finitely presented, and showed that they hold only up to a degree determined by $\# \mathcal{S}$.
These phenomena were extended to arbitrary
reductive groups by the rank theorem of Bux, K\"ohl and Witzel \cite{B}. Important evidence for this was also contributed by
Behr~\cite{Be}, Abels~\cite{Abels}, Abramenko~\cite{AbramenkoFin}, and Bux and Wortman~\cite{Bux-Wortman}.

Let $\mathcal{C}$ be a smooth projective and geometrically integral curve over a finite field $\mathbb{F}$ of characteristic $p$, and let $\mathcal{S}$ be a finite set of closed points of $\mathcal{C}$, of cardinality $s$. Let $k=\mathbb{F}(\mathcal{C})$ be the function field of $\mathcal{C}$, and let $\OS\subset k$ be the ring of functions of $\mathcal{C}$ that are regular outside $\mathcal{S}$. For each $P\in \mathcal{S}$, we write $k_P$ for the completion of $k$ at $P$ and $X_P=X(\mathrm{SL}_2,k_P)$ for the Bruhat-Tits tree of $\mathrm{SL}_2(k_P)$. An $\mathcal{S}$-arithmetic subgroup $G$ of $\mathrm{SL}_2$ is a subgroup of $\mathrm{SL}_2(k)$ that is commensurable with $\mathrm{SL}_2(\OS)$.
Such a group acts diagonally, through simplicial maps, on the product $X=\prod_{P\in \mathcal{S}}X_P,$
which is a contractible simplicial complex of dimension $s.$
The Euler-Poincar\'e characteristic of $G$ is defined by
\begin{equation}\label{eq euler-poincare}
\chi(G)=\sum_{\tau}\dfrac{(-1)^{\mathrm{dim}(\tau)}}{\# \mathrm{Stab_G(\tau)}},
\end{equation}
where $\tau$ runs through a system of representatives of the cells of $\X$ modulo $G$.

Given a family of subgroups $\mathfrak{H}=\{H_\sigma\}_{\sigma \in \Sigma}$ of $G$ and a $G$-module $M$, the relative homology groups $H_i(G \bmod H_{\sigma}, M)$ fit into a long exact sequence
$$
\cdots \to \bigoplus_{\sigma \in \Sigma} H_i(H_\sigma, M) \to H_i(G, M) \to H_i(G \bmod \mathfrak{H}, M) \to \bigoplus_{\sigma \in \Sigma} H_{i-1}(H_\sigma, M) \to \cdots.
$$
Thus they measure how far the homology of $G$ is from being the direct sum of the homologies of the subgroups $H_\sigma$.
A natural family to consider consists of the maximal unipotent subgroups of 
$G$. Indeed, when the torsion of $G$ is $p$-primary, an element of $G$ is unipotent if and only if it has $p$-power order. Hence the maximal unipotent subgroups of $G$ are exactly its maximal $p$-torsion subgroups in this case. Therefore, they are the natural analogues of the $p$-Sylow subgroups of $G$.

Therefore, we are interested in the homology of $G$ relative to a representative system $\{U_\sigma\}_{\sigma\in \Sigma_0}$ of the $G$-conjugacy classes of maximal unipotent subgroups of $G$.
The case $s=1$ is classical and it is addressed by the following result:

\begin{MainThm} \cite[Ch.~II, \S2.9, p.~135 and Thm.~14]{SerreTrees}\label{teoA}
Assume that $s=1$, and $\mathbb{F}$ is finite of characteristic $p$, and let $G$ be a
finite-index subgroup of $\mathrm{SL}_2(\mathcal{O_\mathcal{S}})$ whose torsion is $p$-primary.
Let $\{\xi_\sigma\}_{\sigma\in\Sigma_0}$ be a set of representatives of the
$G$-orbits in the visual limit of $X$, or equivalently the $G$-orbits in
$\mathbb{P}^1(k)$, and write $G_\sigma = \mathrm{Stab}_G(\xi_\sigma)$. Then:
\begin{enumerate}
  \item[(a)] the set $\mathfrak{G} = \{G_\sigma\}_{\sigma\in\Sigma_0}$ is a
  representative system for the $G$-conjugacy classes of maximal unipotent
  subgroups in $G$, and
  \item[(b)] for all $i \neq 1$, the homology groups
  $H_i(G \bmod G_{\sigma}, M)$ are trivial. 
  Moreover, if $M$ is 
  $\mathbb{Z}$-finitely generated, then
  $H_1(G \bmod G_{\sigma}, M) \cong M^{-\chi(G)}$, where $\chi(G)$ is defined as in Eq. \eqref{eq euler-poincare}.
\end{enumerate}
\end{MainThm}

In the context of Th. \ref{teoA}, the complex $X$ is the single tree $X_P$, and Serre computes the relative homology of $G$ by looking at the subcomplex of $X$ consisting of the cells with non-trivial $G$-stabilizer.
The connected components of that subcomplex are in bijection with $\mathbb{P}^1(k)$, and the resulting object is the Steinberg module $\mathrm{St}$, defined as the kernel of the augmentation map $\mathbb{Z}[\mathbb{P}^1(k)] \to \mathbb{Z}$.
For $s=1$ the module $\mathrm{St}$ turns out to be projective, whence the relative homology of $G$ is concentrated in degree one and is completely determined by the Euler-Poincar\'e characteristic of $G$.

In~\cite[\S~3, p.~155]{SerreCohomology}, Serre conjectured that Th. \ref{teoA} extends to any $\mathcal{S}$-arithmetic subgroup of a reductive group $\mathbf{G}_k$ whose $\mathcal{S}$-rank $\mathrm{rk}_{\mathcal{S}}(\mathbf{G}_k)$ is the sum of its local ranks at the places of $S$, or equivalently the dimension of the associated product of Bruhat-Tits buildings. The following statement makes this conjecture precise  (see \cite{BravoLoisel2}).
%In~\cite[\S~3, p.~155]{Se2}, Serre conjectures these results, as well as Theorem~\ref{teo max unip serre}, extend to any $\mathcal{S}$-arithmetic subgroup of a reductive group. Next statement formalizes this Serre's conjecture:

\begin{MainConj}\label{conj serre rel hom}
Let $\mathbf{G}_k$ be a $k$-isotropic reductive group. Let $r=\mathrm{rk}_{\mathcal{S}}(\mathbf{G}_k)$ be the $\mathcal{S}$-rank of $\mathbf{G}_k$.
Let $G$ be a finite-index subgroup of $\mathbf{G}_k(\mathcal{O}_{\mathcal{S}})$ whose torsion is $p$-primary, for $p=\mathrm{char}(\mathbb{F})>0$, and let $\lbrace U_{\sigma} \rbrace_{\sigma}$ be a representative set of the $G$-conjugacy classes of maximal unipotent subgroups in $G$.
Let $M$ be a $G$-module.
Then:
\begin{itemize}
\item[(a)] for all $i > r$, the relative homology groups $H_i (G \bmod U_{\sigma}, M)$ are trivial, and there exist non-positive integers $n_1, \dots, n_r$ such that:
\item[(b)] if $M$ is finitely generated over $\mathbb{Z}$, then, for each $i\in \lbrace 1, \dots, r \rbrace$, we have $H_i(G \bmod U_{\sigma}, M)\cong M^{-n_i}$, and
\item[(c)] $\sum_{i=1}^r (-1)^{i+1} n_i=\chi(G)$, where $\chi(G)$ is the Euler-Poincar\'e characteristic of $G$ defined as in Eq. \ref{eq euler-poincare}.
\end{itemize}
\end{MainConj}

%More recently, the second named author studied in \cite{CB} the relative homology of arithmetic subgroups of $\mathrm{SU}(3)$, and $\mathcal{S}=\lbrace P \rbrace$, following closely the ideas of Serre.
%Since the split rank of $\mathrm{SU}(3)$ is one, the associated building is again a tree, and most of Serre's arguments can be adapted once $\mathbb{P}^1(k)$ is replaced by the appropriate set of rational boundary points (i.e., the $k$-rational points of the Borel variety of $\mathrm{SU}(3)$).

More recently, the second-named author extended Serre's approach in \cite{CB} to the relative homology of arithmetic subgroups of $\mathrm{SU}(3)$, with $\mathcal{S}=\lbrace P \rbrace$, replacing $\mathbb{P}^1(k)$ by the $k$-rational points of the Borel variety of $\mathrm{SU}(3)$.
That work generalizes Serre's results and gives a positive answer to Conj. \ref{conj serre rel hom}, covering all the semisimple simply connected quasi-split algebraic groups with split rank $1$.
In the present article we take the orthogonal route of generalization. We keep $\mathrm{SL}_2$, but we allow $\mathcal{S}$ to have arbitrary cardinality, so that the associated building $\X$ is the product of the $s$ trees introduced above.

Even for $s=1$, relative homology has important applications. Teitelbaum used in \cite{teitelbaum} to understand suitable spaces of modular forms of Drinfeld type.
More precisely, Teitelbaum constructs an isomorphism between a space of Drinfeld cuspidal forms and a space of harmonic cocycles. 

\section{Main results}\label{Section main results}

%Let us now make precise the object we compute.
In order to state our main results precisely, we introduce some notations.
Let $R$ be a commutative ring with unity and let $\mathfrak{H}=\lbrace H_{\sigma} \rbrace_{\sigma \in \Sigma}$ be a non-empty family of subgroups of $G$.
Let $\varepsilon : \bigoplus_{\sigma \in \Sigma} R[G/H_{\sigma}] \to R$ be the augmentation map and set $\Delta_{\mathfrak{H}}=\ker(\varepsilon)$.
For each $R[G]$-module $M$ and each $i \in \mathbb{Z}_{\geq 1}$, the $M$-valued relative homology groups of $G$ modulo $\mathfrak{H}=\lbrace H_{\sigma} \rbrace_{\sigma \in \Sigma}$ are defined by
\[H_i\bigl(G \bmod H_{\sigma},M\bigr):=
\mathrm{Tor}_{i-1}^{R[G]}(\Delta_{\mathfrak{H}},M) \cong
H_{i-1}\bigl(G, \Delta_{\mathfrak{H}} \otimes_{R} M\bigr),\]
and the relative $R$-homological dimension is defined by $\mathrm{rhd}_R(G,H_{\sigma}):=\mathrm{fd}_{R[G]}(\Delta_{\mathfrak{H}})+1$, where $\mathrm{fd}$ denotes the flat dimension of $\Delta$ as an $R[G]$-module. 
In particular, $H_r(G \bmod H_\sigma, M) = 0$ for every $R[G]$-module $M$ and every integer $r > \mathrm{rhd}_R(G, H_\sigma)$, so that $H_R(G,M)\cong \bigoplus_{\sigma \in \Sigma_0} H_r(H_{\sigma},M)$, for all $r> \mathrm{rhd}_R(G, H_\sigma)$.

Let $\Sigma_0$ be a set of representatives of the $G$-orbits in
$\mathbb P^1(k)$. For each $\sigma\in\Sigma_0$, let
$G_\sigma=\operatorname{Stab}_G(\sigma)$, and let $U_\sigma$ be the
subgroup of unipotent elements of $G_\sigma$. The family
$\mathfrak U=\{U_\sigma\}_{\sigma\in\Sigma_0}$ is a representative system of the
$G$-conjugacy classes of maximal unipotent subgroups of $G$
\cite[Thm.~3.1]{BravoLoisel2}. Write
$\mathfrak G=\{G_\sigma\}_{\sigma\in\Sigma_0}$. The orbit maps
$gG_\sigma\mapsto g\sigma$ induce an augmentation-preserving isomorphism
\[
  \bigoplus_{\sigma\in\Sigma_0} R[G/G_\sigma]
  \xrightarrow{\ \sim\ } R[\mathbb P^1(k)].
\]
The kernel of $R[\mathbb P^1(k)]\to R$ is the Steinberg module $\mathrm{St}=\mathrm{St}_R$. Thus $\Delta_{\mathfrak G}$ identified with $\mathrm{St}$.
For the unipotent family, the corresponding module is
$\Delta_{\mathfrak U}=\ker(\bigoplus_{\sigma\in\Sigma_0}R[G/U_\sigma]\to R)$.
Consequently, for every $R[G]$-module $M$ and every $i\ge1$,
\[
\begin{aligned}
 H_i(G\bmod G_\sigma,M)&\cong
   \operatorname{Tor}^{R[G]}_{i-1}(\mathrm{St}_R,M),\\
 H_i(G\bmod U_\sigma,M)&\cong
   \operatorname{Tor}^{R[G]}_{i-1}(\Delta_{\mathfrak U},M).
\end{aligned}
\]
Here both augmentation kernels are viewed as right $R[G]$-modules via
inversion.

Unless otherwise stated, we assume that the torsion of $G$ is $p$-primary.
Our first main result exposes a discrepancy with the case $s=1$, in which the Steinberg module is projective.

\begin{theorem}\label{main theorem intro pd}
The projective and the flat dimension of the Steinberg module $\mathrm{St}$ over $R[G]$ are both equal to $s-1$.
In particular $\mathrm{rhd}_R(G,G_{\sigma})=s$.
\end{theorem}

An analogous result follows for $\Delta_{\mathfrak{U}}$.

\begin{theorem}\label{main theo intro pd of delta}
If $p\in R^{*}$, then the flat dimension of $\Delta_{\mathfrak{U}}$ over $R[G]$ is $s-1$.
In particular $\mathrm{rhd}_R(G,U_{\sigma})=s$.
\end{theorem}

%Rephrased in terms of the relative homology, Theorem \ref{theorem intro pd} says that \[ \mathrm{rhd}_R(G,G_{\sigma})=s, \]
Note that Th. \ref{main theorem intro pd} implies that the relative homology of $G$ modulo $G_{\sigma}$ may be non-trivial in every degree $1 \leq i \leq s$, and vanishes above $s$.

The remaining possibly non-trivial homological groups can be described in terms of the Euler-Poincar\'e characteristic defined from the action of $G$ on $\X$. 
Let us write $\Xinf$ for the subcomplex of $\X$ consisting of the cells with non-trivial $G$-stabilizer, and, for each $i$, let $S_i$ be the set of $i$-dimensional cells of $\X$ outside $\Xinf$.
Since $G$ acts freely on $S_i$, the $R[G]$-module $L_i=R[S_i]$ is free, of rank $l_i=\#(G \backslash S_i)$. See \S \ref{section geometry of Xinf} for details.
Our second main result expresses $\chi(G)$ purely in terms of these ranks.

\begin{theorem}\label{main theorem intro chi}
The cells with non-trivial stabilizer contribute nothing to $\chi(G)$. Moreover, for each $i$ the rank $l_i$ is finite, and one has the identity
\[\chi(G)= \sum_{i=0}^{s} (-1)^{i}l_i.\]
\end{theorem}

%The content of Theorem \ref{theorem intro chi} is that the cells with non-trivial stabilizer contribute nothing to $\chi(G)$.
Similar results can be found in \cite{SerreTrees,CB} in the context of $s=1$.
In all cases, the cells with non-trivial stabilizers organize into sectors along which the stabilizer controlled increases, so that the contribution faces in a given alcove cancel out (in pairs).
Combining the two previous results, we obtain the description of relative homology in terms of $\chi(G)$, which is our third main result.

\begin{theorem}\label{main theorem intro euler}
Let $D$ be a principal ideal domain and let $M$ be a $D[G]$-module whose underlying $D$-module has finite rank $m:=\rk_{D}(M)<\infty$. Then:
\begin{itemize}
  \item[(a)] $H_i(G\bmod G_\sigma,M)$ has finite rank for every $i$, and vanishes for
        $i>s$ and for $i=0$, and
  \item[(b)] one has the identity
        $$\sum_{i=1}^{s}(-1)^{i} \rk_{D} H_i\big(G\bmod G_\sigma,M\big) =\chi(G)\cdot m .$$
\end{itemize}
\end{theorem}

Using Th. \ref{main theorem intro euler}, we describe the relative homology groups $H_r(G \bmod U_{\sigma},M)$, for $r\leq s$, in terms of $\chi$ and the module of coinvariants $M_H=H_0(H,M)=M / \langle hm-m \vert h\in H, m\in M \rangle$.

\begin{theorem}\label{main thm euler-unipotent for U}
Assume that $s \geq 2$ and that the torsion of $G$ is $p$-primary. Let $D$ be
a principal ideal domain such that $p \in D^{*}$, and let $M$ be a $D[G]$-module of
finite rank $m := \mathrm{rk}_D(M) < \infty$. Then:
\begin{enumerate}
\item[(a)] $H_i(G \bmod U_\sigma, M)$ has finite rank for every $i$, and
vanishes for $i > s$.
\item[(b)] One has the identity
\[
\sum_{i=1}^{s} (-1)^i \mathrm{rk}_D H_i(G \bmod U_\sigma, M)
= \chi(G) \cdot m - \sum_{\sigma \in \Sigma_0} \mathrm{rk}_D\big(M_{U_\sigma}\big).
\]
\end{enumerate}
\end{theorem}

Note that when $M=D$ with the trivial $G$-action, the preceding result shows that the missing term between the alternating sum of ranks of $H_r(G \bmod U_{\sigma},M)$, for $r \leq s$, and $\chi(G)$ is the number $h$ of cusps of $G$, which is closely related to the class number $\# \mathrm{Pic}(\OS)$ (see Corollary~\ref{coro euler-unipotent-trivial}).
In particular, Conjecture~B~(c) does not hold as stated for the maximal unipotent subgroups when $s \geq 2$.

The flat dimension of $R$ as an $R[G]$-module with trivial $G$-action is called the homological dimension of $G$ and denoted by $\mathrm{hd}_R(G)$ (cf. \cite[\S 4.4]{berihomologicaldimensiondisgroups}).
Equivalently, $\mathrm{hd}_R(G)$ is the smallest integer $n \geq 0$ such that, for each $R[G]$-module $M$, we have $H_r(G,M) = \lbrace 0 \rbrace$ for all $r>n$.
In the number field case, Borel and Serre computed the (co)homological dimension of $S$-arithmetic groups in \cite{Borel-Serre}. In positive characteristic the situation is more subtle, since $S$-arithmetic groups always contain $p$-torsion, so one only considers coefficient rings $R$ with $p \in R^{*}$.
In that setting, the action on the Bruhat-Tits building $X_{\mathcal{S}}$ gives the $\mathcal{S}$-rank $\mathrm{rk}_S(\mathbf{G}_k)$ as upper bound \cite{SerreCohomology}, which Hughes showed to be optimal over $\mathbb{Q}$ for split semisimple simply connected groups \cite{Hughes2023}.
However, the case of arbitrary coefficient modules, such as $\mathbb{Z}/l\mathbb{Z}$, with $l \neq p$, remains largely open.
The next result describes the homological dimension of $S$-arithmetic subgroups of $\mathrm{SL}_2$ showing that the upper bound is optimal in this case even if $R=\mathbb{Z}/l\mathbb{Z}$ with $l\neq p$.

\begin{theorem}\label{main teo HD}
Let $G$ be an $\mathcal{S}$-arithmetic subgroup of $\mathrm{SL}_2(k)$. Assume that the torsion of $G$ is $p$-primary and $p\in R^{*}$. Then $\mathrm{hd}_R(G)=s$
\end{theorem}

Our approach relies  on the subcomplex $\Xinf \subseteq \X$ of cells with non-trivial
stabilizer. The stabilizer of each vertex $v \in \Xinf$ is a non-trivial finite
$p$-group, so it fixes a unique $D_v \in \mathbb{P}^1(k)$ and fixes pointwise the
sector chamber $Q(v, D_v)$. This identifies the connected components of $\Xinf$
with $\mathbb{P}^1(k)$, whence $H_1(\X \bmod \Xinf; R) \cong \mathrm{St}$
(Prop.~\ref{lemma H_1}). Moreover, any finite family of sector chambers with the same
direction meets in a sector chamber, so the nerve lemma shows that each component
is contractible. As a consequence, the complex $L_\bullet$ of relative chains of
$(\X, \Xinf)$, which consists of finitely generated free $R[G]$-modules, yields a
projective resolution of $\mathrm{St}$ of length $s-1$. This gives the upper bound in
Theorem~\ref{main theorem intro pd} and identifies $H_*(G \bmod G_\sigma, M)$ with the
$G$-equivariant homology of $(\X, \Xinf)$. Theorem~\ref{main theorem intro euler} then
follows by computing the Euler characteristic of $L_\bullet \otimes_{D[G]} M$ in two
ways.
By using that and Lyndon-Holchild-Serre spectral sequence, we prove Ths. \ref{main theo intro pd of delta} and \ref{main thm euler-unipotent for U} in \S \ref{section Rel. Homology for U}.

We close this section by describing two directions in which the present work could be contributed.
Nothing above is specific to $\mathrm{SL}_2$ beyond the fact that its rational spherical building is $\mathbb{P}^1(k)$, so it should be possible to merge the present article with \cite{CB} and to treat $\mathcal{S}$-arithmetic subgroups of $\mathrm{SU}(3)$ for an arbitrary $\mathcal{S}$, extending our results to all semisimple simply connected quasi-split algebraic groups of rank $1$.

On the other hand, more speculatively, for $s=1$ the relative homology is Teitelbaum's bridge between chains and Drinfeld cuspidal modular forms, whereas for $s \geq 2$ there is, to our knowledge, no theory of Drinfeld modular forms associated with the action of $G$. Nevertheless, it would be interesting to consider a definition of Drinfeld modular forms for the action of $G$ and to see whether the space of such forms can be related to automorphic forms from an adelic perspective via relative homology.
We plan to pursue 
this idea, along with its consequences, in future works.

\subsection*{Overview of the article}
Section \ref{section context and main results} fixes the notation, collects the geometric properties of the product of trees $\X$ that we shall need, in particular the behavior of sector chambers and the action of finite $p$-subgroups on $\mathbb{P}^1(k)$, and recalls the definitions of the Steinberg module and of the relative homology of a group modulo a family of subgroups.
Section \ref{section geometry of Xinf} is devoted to the subcomplex $\Xinf$. There we describe its connected components in terms of $\mathbb{P}^1(k)$ and we prove, via the nerve lemma, that each of them is contractible.
Section \ref{section res and dim st} constructs the resolution of $\mathrm{St}$ by relative chains and proves Theorem \ref{main theorem intro pd}.
Section \ref{section Rel. Homology} computes the relative homological dimension of $G$, identifies the relative homology of $G$ modulo its cusps with the $G$-equivariant homology of the pair $(\X,\Xinf)$, and proves Theorems \ref{main theorem intro chi} and \ref{main theorem intro euler}.
Finally, Section~\ref{section Rel. Homology for U} compares the relative homology of $G$ modulo the stabilizers $G_\sigma$ with that modulo the maximal unipotent subgroups $U_\sigma$, proves that $\mathrm{fd}_{R[G]}(\Delta_{\mathfrak{U}})$ is $ s-1$ when $p \in R^{*}$, and computes the alternating sum of the ranks of $H_i(G \bmod U_\sigma, M)$ in terms of $\chi(G)$ and the coinvariants $M_{U_\sigma}$.

\begin{comment}
Let $\mathbf{G}$ be a \textbf{semi-simple} and simply connected Chevalley group (scheme) defined over $\mathbb{Z}$. Let $\mathcal{C}$ be a smooth, projective, geometrically integral curve over a \textbf{finite} field $\mathbb{F}$. Let $\wp$ be a closed point on $\mathcal{C}$, and let us denote by $A$ the ring of functions that are regular outside $\lbrace \wp \rbrace$. We denote by $k$ the fraction field of $A$.

Let us denote by $\mathbf{G}(A)$ the abstract group consisting in the $A$-points of $\mathbf{G}$.
This acts on the euclidean Bruhat-Tits building $X$ defined from $\mathbf{G}$ and the completion $K$ of $k$ at $\wp$. Let $G$ be a finite index subgroup of $\mathbf{G}(A)$ and assume that $G$ is \textit{without $p'$-torsion}, i.e. where the elements of $G$ of finite order have order a power of $p$.

By following \cite[??]{??} we define the Euler-Poincar\'e characteristic $\chi(G)$ of $G$ by the identity
\begin{equation}\label{eq euler-poincare char}
\chi(G)=  \sum_{F }  \frac{(-1)^{\mathrm{dim}(F)}}{ \mathrm{Card}(\mathrm{Stab}_G(F))},
\end{equation}
where $F$ runs through a system of representative of the cellular complex of $X$ modulo $G$, i.e. $\sigma$ runs through a suitable lift of $G \backslash X$. It follows from \cite[Prop. 13.8]{B} that $\chi(G)$ is finite. In the present work, we interpret $\chi(G)$ in terms of the action of $G$ on the spherical building $\partial_{\infty}(X)$. Then, we use this in order to describe the relative homology of $G$ relative to a certain family of subgroups in terms of $\chi(G)$. 
\end{comment}

\section{Preliminaries and the geometry of Bruhat-Tits buildings}\label{section context and main results}

Let $\mathcal{C}$ be a smooth, projective and geometrically integral curve over a finite field $\mathbb{F}$, and let $\mathcal{S}$ be an arbitrary finite set of closed points in $\mathcal{C}$ of cardinality $s$. Let $k=\mathbb{F}(\mathcal{C})$ be the function field of $\mathcal{C}$ and let $\mathcal{O}_{\mathcal{S}}\subset k$ be the ring of functions that are regular outside $\mathcal{S}$. Throughout this article we assume that $\mathbb{F}$ is finite of characteristic $p> 0$.

We say that $G \subset \mathrm{SL}_2(k)$ is an $\mathcal{S}$-arithmetic subgroup of $\mathrm{SL}_2(k)$ when $G$ is commensurable with $\mathrm{SL}_2(\mathcal{O}_\mathcal{S})$.
Moreover, we say that the torsion of $G$ is $p$-primary whenever any torsion element in $G$ has a $p$-power order.

\subsection{On the geometry of \texorpdfstring{$X$}{X}}\label{section visual limit}

For each $P \in \mathcal{S}$, we denote by $k_P$ the completion of $k$ at $P$.
The valuation of $k$ defined by $P$ extends to a valuation $\nu_P:k_P\to \mathbb{Z}\cup \lbrace \infty \rbrace$.
Let us fix a uniformizing parameter $\pi=\pi_P$ which is a generator of the unique maximal ideal of $\mathcal{O}_P=\lbrace x \in k_P \vert \nu_P(x) \geq 0 \rbrace$.
The Bruhat-Tits tree $X_P=X(\mathrm{SL_2},k_P)$ of $\mathrm{SL}_2$ at $k_P$ is the graph defined as follows:
\begin{itemize}
  \item vertices are homothety classes $[L]$ of $\mathcal{O}_P$-lattices $L \subset k_P^2$, i.e., free $\mathcal{O}$-submodules of rank $2$, where $L \sim cL$ for $c \in k_P^\times$;
  \item edges join $[L]$ and $[L']$ whenever they admit representatives with $\pi_P L \subsetneq L' \subsetneq L$.
\end{itemize}
The tree $X_P$ is a $(q+1)$-regular tree, where $q= \#(\mathcal{O}_P/\pi_P\mathcal{O}_P)$. Indeed, the set of neighbors of a vertex is in bijection with $\mathbb{P}^1(\mathcal{O}_P/\pi_P\mathcal{O}_P)$.
The group $\mathrm{SL}_2(k_P)$ acts on $X_P$ by
isometries without inversions of edges.

The visual limit $\partial(R)$ of a ray $R$ in $X_P$ is the class of all rays parallel to $R$.
The set of visual limits of all rays in $X_P$ is in bijection with $\mathbb{P}^1(k_P)$.
The action of $\mathrm{SL}_2(k_P)$ in $X_P$ induces an action of $G$ on $\mathbb{P}^1(k_P)$, which is exactly the action by M\"{o}bius transformations.

Any $\mathcal{S}$-arithmetic group $G \subset \mathrm{SL}_2(k)$ acts diagonally on the $s$-dimensional Bruhat-Tits building $X:=\prod_{P\in \mathcal{S}}X_P$.
The Bruhat-Tits building $X$ is a CW-complex.
An alcove of $X$ is the product of exactly $s$ edges, each one in the corresponding tree $X_P$. A facet of $X$ is a face of some alcove. The next result describes the stabilizer of $G$ acting on $X$.

\begin{lemma}\label{lemma estabilizadores finitos}
Let $v$ (or $\sigma$) be a vertex (or a facet) of $X$, then $\mathrm{Stab}_G(v)$ is finite. In particular, the stabilizers are torsion groups.
\end{lemma}
\begin{proof}
Recall that $G$ is commensurable with $\mathrm{SL}_2(\OS)$.
Hence $\mathrm{Stab}_G(v)$ is commensurable with $\mathrm{Stab}_{\mathrm{SL}_2(\OS)}(v)$, whence we can assume that $G=\mathrm{SL}_2(\OS)$. 
The result then follows from \cite[Obs. 10.1]{B}.
\end{proof}

A sector chamber $Q$ in $\X$ is the product $Q=\prod_{P \in \mathcal{S}} R_P$, where $R_P$ is a ray in $X_P$.
The visual limit of $Q$ is $\prod_{P \in \mathcal{S}} \partial(R_P) \in \prod_{P \in \mathcal{S}} \mathbb{P}^1(k_P)$.

\begin{lemma}\label{lemma equal visual limit}
Let $Q$ and $Q'$ be two sector chambers whose visual limits coincide. Then $Q \cap Q'$ is a sector chamber.
\end{lemma}

\begin{proof}
We can write $Q=\prod_{P \in \mathcal{S}} Q_P$ and $Q'=\prod_{P \in \mathcal{S}} Q'_P$, where $Q_P$ and $Q_P'$ are two rays in $X(\mathrm{SL}_2, k_P)$ whose visual limits coincide.
Then $Q_P \cap Q_P'$ is a ray, say $Q''_P$, in $X(\mathrm{SL}_2, k_P)$.
Therefore $Q \cap Q'=\prod_{P \in \mathcal{S}} Q_P''$ is a sector chamber of $\X$.
\end{proof}

\begin{corollary}\label{cor:finite-int}
Let $D\in\mathbb P^1(k)$ and let $Q_1, \dots, Q_n$ be finitely many sector chambers whose visual limit is $D$. Then $\bigcap_{i=1}^{n} Q_{i}$ is a sector chamber of $\X$.
\end{corollary}

Note that the group $G$ acts on $\mathbb{P}^{1}(k)$ by M\"{o}bius transformations. 
This action is compatible with the action induced by the diagonal inclusion $\mathbb{P}^1(k) \to \prod_{P \in \mathcal{S}} \mathbb{P}^1(k_P)$.
The next result is well known and follows from \cite[Ch. II, \S 2.9]{SerreTrees}. The argument below is entirely standard, but we include it only for completeness.

\begin{lemma}\label{Lemma action of a p-group}
Let $P$ be a non-trivial finite $p$-subgroup of $\mathrm{SL}_2(k)$.
Then the natural action of $P$ on $\mathbb{P}^1(k)$ has a unique fixed point.
\end{lemma}

\begin{proof}
Let $P$ be a non-trivial finite $p$-subgroup of $\mathrm{SL}_2(k)$.
Since $\mathscr{G}(k)/\mathscr{B}(k) \cong \mathbb{P}^1(k)$,
we must prove that there exists a unique non-zero vector in $k^2$ fixed by $P$. We first establish the existence of such a vector. Let $v \in k^{2}$ be a non-zero vector. Define $W$ as the $\mathbb{F}$-vector space generated by all elements in the $P$-orbit of $v$.
One checks at once that $P \cdot W= W$.
Since $P$ is finite, the orbit $P\cdot v$ is finite as well, whence $W$ is finite dimensional. In particular, $\# W=p^n$, for some $n \in \mathbb{Z}_{>0}.$
Let us write $W^P:=\lbrace x \in W : P \cdot x = \lbrace x \rbrace \rbrace$.
Since $0 \in W^{P}$, we have $1 \leq \# W^P \leq p^n$.
Moreover, since $W^P$ is exactly the subset of $W$ where $P$ acts trivially, an orbit count yields
\begin{equation}\label{eq W}
    \# W = \# W^P + \sum_{i=1}^{m} \# \mathrm{Orb}_{P}(x_i),
\end{equation}
where $\# \mathrm{Orb}_{P}(x_i) >1$.
Note that each $\# \mathrm{Orb}_{P}(x_i)=\frac{\# P}{ \# \mathrm{Stab}_P(x_i)}$ is divisible by $p$.
Thus, it follows from Eq. \eqref{eq W} that $\# W^P$ is also divisible by $p$.
We conclude that there exists a non-zero $P$-fixed point in $W \subseteq k^{2}$.
We now prove uniqueness. To this end, assume that $P$ fixes two points $\xi_1, \xi_2 \in \mathbb{P}^1(k)$.
In other words, assume that
\begin{equation}\label{eq P subset stab}
P \subseteq \mathrm{Stab}_{\mathrm{SL}_2(k)}(\xi_1) \cap \mathrm{Stab}_{\mathrm{SL}_2(k)}(\xi_2).
\end{equation}
Since $\mathrm{SL}_2(k)$ acts transitively on $\mathbb{P}^1(k)$, we can write $\xi_1= g_1 \cdot{\infty}$ and $\xi_2= g_2 \cdot {\infty}$, for some $g_1, g_2 \in \mathrm{SL}_2(k)$.
By the Bruhat decomposition, we can also write $g_1^{-1}g_2= b \mathrm{s}^{\epsilon} b'$, where $b,b' \in \mathrm{B}(k)=\ssbmatrix{k}{k}{0}{k} \cap \mathrm{SL}_2(k)$, and $\epsilon \in \lbrace 0,1 \rbrace$.
Assume first that $\epsilon=0$. Then $g_1^{-1}g_2\in \mathrm{B}(k)$, which equals $\mathrm{Stab}_{\mathcal{G}(k)}(\infty)$. Hence $\xi_1= g_1 \cdot {\infty}= g_1 (g_1^{-1}g_2) \cdot {\infty}= g_2 \cdot {\infty}= \xi_2$.
Assume now that $\epsilon=1$. It follows from Identity~\eqref{eq P subset stab} that $P \subseteq g_1\mathrm{B}(k) g_1^{-1} \cap g_2  \mathrm{B}(k) g_2^{-1}$. Since $g_1^{-1}g_2= b \mathrm{s} b'$, with $b,b' \in \mathrm{B}(k)$, we have
$$(g_1b)^{-1} P (g_1b) \subseteq \mathrm{B}(k) \cap \mathrm{s}  \mathrm{B}(k) \mathrm{s}^{-1}.$$
Therefore, since $\mathrm{B}(k) \cap \mathrm{s}  \mathrm{B}(k) \mathrm{s}^{-1}=\mathrm{T}(k)=\ssbmatrix{k}{0}{0}{k} \cap \mathrm{SL}_2(k)$, we obtain a $p$-subgroup contained in the diagonal torus $\mathrm{T}(k)\cong k^{*}$, which is impossible.
\end{proof}

The Euler-Poincar\'e characteristic $\chi(G)$ of $G$ is defined as follows
\begin{equation}\label{eq euler-poincare char}
\chi(G)=  \sum_{\sigma }  \frac{(-1)^{\mathrm{dim}(\sigma)}}{ \# \mathrm{Stab}_G(\sigma)},
\end{equation}
where $\sigma$ runs through a system of representatives of the cells of $\X$ modulo $G$, i.e., through a suitable lift of $G \backslash \X$.
It follows from \cite[Prop. 13.8]{B} that $\chi(G)$ is finite.

\subsection{The Steinberg module and the relative homology}\label{section Stainbert}

Let $R$ be a commutative unital ring and let $R[\mathbb{P}^1(k)]$ be the free $R$-module with basis $\mathbb{P}^1(k)$.
Let $\epsilon: R\, [\mathbb{P}^1(k)] \to R$ be the augmentation homomorphism determined by $\epsilon(z)=1$ for every $z \in \mathbb{P}^1(k)$.
We define the Steinberg module $\mathrm{St}=\St_R$ associated to the action of $G$ on $\mathbb{P}^1(k)$ as the kernel of $\epsilon$.
Thus, we have the exact sequence
\begin{equation}\label{Eq ex sequence St}
0 \to \mathrm{St} \to R[\mathbb{P}^1(k)] \to R \to 0.
\end{equation}
Since $\mathrm{SL}_2(k)$ acts on $\mathrm{St}$, the latter is an $R[G]$-module.

We now recall the definition of the relative homology of a group $G$ with respect to a non-empty family of subgroups $\mathfrak{H}=\lbrace H_{\sigma} \rbrace_{\sigma \in \Sigma}$.

\begin{definition}\label{def relative homology}
Let $R$ be a commutative unital ring, $G$ a group and $\mathfrak{H}=\lbrace H_{\sigma}\rbrace_{\sigma\in\Sigma}$
a non-empty family of subgroups of $G$.
Let $\varepsilon: \bigoplus_{\sigma\in\Sigma} R[G/H_{\sigma}]\to R$ be the augmentation map defined by $gG_{\sigma}\mapsto 1$. Set $\Delta=\Delta_{\mathfrak{H}}=\ker(\varepsilon)$,
considered as a right $R[G]$-module via $g\mapsto g^{-1}$.

For each $R[G]$-module $M$ and each $i\in\mathbb{Z}_{\geq 1}$, the $M$-valued relative homology groups of $G$ modulo
$\mathfrak{H}=\lbrace H_{\sigma}\rbrace_{\sigma\in\Sigma}$ are defined by
$$H_i\bigl(G \text{ mod } H_{\sigma},M\bigr):=
\mathrm{Tor}_{i-1}^{R[G]}(\Delta,M) \cong
H_{i-1}\bigl(G, \Delta \otimes_{R} M\bigr),$$
where $\Delta\otimes_R M$ is endowed with the diagonal $G$-action. Moreover, we set $H_0(G \text{ mod } H_{\sigma},M):=0$.

The relative $R$-homological dimension of $G$ modulo
$\mathfrak{H}=\lbrace H_{\sigma}\rbrace_{\sigma\in\Sigma}$ is
$$\mathrm{rhd}_R(G,H_{\sigma}):=\mathrm{fd}_{R[G]}(\Delta)+1
\ \in\ \mathbb{Z}_{\geq 0}\cup\lbrace\infty\rbrace,$$
the flat dimension of $\Delta$ as an $R[G]$-module, shifted by one.
Equivalently, $\mathrm{rhd}_R(G,H_{\sigma})$ is the smallest $n\geq 0$ such that
$H_r(G\text{ mod }H_{\sigma},M)=\lbrace 0\rbrace$ for every $R[G]$-module $M$
and every $r>n$ (and $\infty$ if no such $n$ exists).
\end{definition}

In the sequel, we assume that $G$ is an $\mathcal{S}$-arithmetic subgroup of $\mathrm{SL}_2(k)$ and that $\mathfrak{G}=\lbrace G_{\sigma} \rbrace_{\sigma \in \Sigma}$ is the set of $G$-stabilizers of a representative system of $G \backslash \mathbb{P}^1(k)$

\begin{lemma}\label{lemma Deleta=St}
With the above notation, $\Delta_{\mathfrak{G}}$ is the Steinberg module $\mathrm{St}$ defined in Equation~\eqref{Eq ex sequence St}. 
\end{lemma}

\begin{proof}
By the orbit-stabilizer relation we have $R[\mathbb{P}^1(k)] = R\big[ \coprod_{\sigma \in \Sigma_0} G/G_{\sigma} \big]= \bigoplus_{\sigma \in \Sigma_0} R[G/G_{\sigma}],$
as an $R[G]$-modules.
Hence, the module $\Delta=\mathrm{ker}(\varepsilon)$ is none other than the Steinberg module $\mathrm{St}$. 
\end{proof}

\section{On the geometry of the subcomplex \texorpdfstring{$\Xinf
\subseteq \X$}{Xi}}\label{section geometry of Xinf}

In this section, we introduce a suitable subcomplex of $\X$ that allows us to understand the Steinberg module defined in \S \ref{section context and main results}.
This enables us to determine the relative homology groups of $G$ modulo the cusp stabilizers defined from the action of $G$ on $\X$ and then modulo a representative system of the $G$-conjugacy classes of its maximal unipotent subgroups.

\begin{definition}
We say that a simplex $\sigma \subset \X$ is $G$-stable if its $G$-stabilizer is trivial, i.e., $\mathrm{Stab}_G(\sigma)=\lbrace e \rbrace$.
In any other case we say that $\sigma$ is $G$-unstable.
Note that if $\sigma$ is a $G$-stable simplex, then any other simplex containing $\sigma$ is also $G$-stable.
Hence the $G$-unstable cells form a subcomplex, which we denote by $\Xinf$.
For each $i \in \lbrace 0, \dots, s \rbrace,$ we write
\[S_i= \big\lbrace  \sigma \text{ simplex }: \mathrm{dim}(\sigma)=i,\, \sigma \subset \X \setminus \Xinf \big\rbrace.\]
By definition of $\Xinf$, the group $G$ acts freely on each $S_i$.
For each index $i$, also write $L_i$ for the $R[G]$-module $L_i= R[S_i]$, whose dimension is exactly $l_i= \#(G \backslash S_i)$.
\end{definition}

\begin{proposition}\label{prop li finite}
Let $G$ be a subgroup of $\mathrm{SL}_2(k)$ commensurable with
$\Gamma=\mathrm{SL}_2(\mathcal O_{\mathcal S})$. For $0\le i\le s$,
\[
 \sum_{[F]\in G\backslash X^{(i)}_{\mathrm{cell}}}\frac1{|G_F|}<\infty.
\]
In particular, for every integer $N\ge1$ there are only finitely many
$G$-orbits of $i$-cells $F$ with $|G_F|\le N$. Hence
$l_i=|G\backslash S_i|$ is finite.
\end{proposition}
\begin{proof}
Let $\mathscr C$ be the set of chambers of $X$. By
\cite[Obs.~13.2 and the proof of Prop.~13.8, Eq.~(13)]{B}, the sum of
the reciprocal chamber-stabilizer orders converges for the image of
$\Gamma$ in $\operatorname{Aut}(X)$. The action has finite kernel, and
passing from the image to $\Gamma$ divides each summand by the order
of that kernel.

For a group $J$ acting on $\mathscr C$ with finite stabilizers, put
$m_J=\sum_{[C]\in J\backslash\mathscr C}|J_C|^{-1}$. If $H\le J$ has
finite index, then $m_H=[J:H]m_J$. Indeed, the $H$-orbits in the
$J$-orbit of $C$ correspond to the $J_C$-orbits on $H\backslash J$.
The stabilizer of $Hj$ in $J_C$ is conjugate to $H_{jC}$. Dividing the
orbit-counting identity for this finite set by $|J_C|$ gives the
asserted formula on each $J$-orbit of chambers. Applying it twice with
$H=G\cap\Gamma$ gives $m_G<\infty$.

Now count incident pairs $(F,C)$ with $F$ an $i$-face of $C$, weighting
each $G$-orbit by the reciprocal of its stabilizer order. A chamber is
an $s$-cube and has $b_i=2^{s-i}\binom{s}{i}$ faces of dimension $i$.
If $n(F)$ is the number of chambers containing $F$, orbit-stabilizer
applied first to the faces of a chamber and then to the chambers
containing a face gives
\[
 b_i m_G
 =\sum_{[(F,C)]}\frac1{|G_F\cap G_C|}
 =\sum_{[F]\in G\backslash X^{(i)}_{\mathrm{cell}}}
       \frac{n(F)}{|G_F|}.
\]
All terms are non-negative, and $n(F)\ge1$ for every $F$. This proves
convergence of the stated series. Each orbit with $|G_F|\le N$
contributes at least $1/N$, which proves both finiteness assertions.
\end{proof}

We now interpret each $R[G]$-module $L_i$ in terms of chains in $\X$ relative to $\Xinf$.
More precisely, $L_i$ is naturally identified with the group of $i$-chains of the pair $(\X, X_\infty)$ with coefficients in $R$, that is, the group of $i$-chains in $\X$ modulo those in $\Xinf$, with coefficients in $R$.
Consequently, the long exact sequence in (topological) relative homology of the pair $(\X,\Xinf)$ yields the following long exact sequence, with top dimension $s$:
\begin{equation}\label{eq top rel homology}
\xymatrix{  0 \ar[r]& H_{s}\big(X_{\infty}; R\big)\ar[r]& 
 H_{s}\big(\X ; R \big)\ar[r]&  H_{s}\big(\X \text{ mod } \Xinf ; R \big) \ar[r]& \cdots
\ar`r[d]`[llll]`[lllld]`[lllldr] [dlll]  \\
{} & H_0\big( X_\infty ; R \big) \ar[r]& H_{0}\big(\X ;R \big)\ar[r]&  H_{0}\big(\X \text{ mod } \Xinf ; R \big) \ar[r]& 0 .}
\end{equation}
Since $\X$ is path-connected and contractible, we have $H_i(\X ;R)= \lbrace 0 \rbrace$, for all $i >0$, and $H_0\big(\X ;R \big)= R$. Hence we deduce from the preceding long exact sequence that:
\begin{itemize}\label{eq homology of X}
    \item[(I)] $H_{s}\big( \Xinf ;R \big)=0$,
    \item[(II)] $H_{i}\big(\X \text{ mod } \Xinf  ;R \big) \cong H_{i-1}\big(\Xinf ;R \big)$, for all $i \in \lbrace 2, \dots, s \rbrace$, and
    \item[(III)] $0 \to H_1\big(\X \text{ mod } \Xinf ; R \big) \to H_0\big( \Xinf ;R \big) \to R \to 0$ is exact.
\end{itemize}

In particular, it follows from the Statement (III) that $H_1(\X \text{ mod }\Xinf ;R)$ may be identified with the kernel of the augmentation homomorphism $H_0(\Xinf ;R) \to R$. The following lemma proves that this kernel is, in fact, the Steinberg module defined by Equation \eqref{Eq ex sequence St}.

\begin{proposition}\label{lemma H_1}
The $0$-homology group $H_0\big( \Xinf ;R\big)$ is isomorphic to $R[\mathbb{P}^1(k)]$. In particular, we have $H_1\big( \X \text{ mod }\Xinf ;R \big) \cong \mathrm{St}$.
\end{proposition}

In order to prove this proposition, we first introduce some notation. Let $v$ be a vertex of $\Xinf$. Then, by definition of $\Xinf$, we have $\mathrm{Stab}_G(v) \neq \lbrace e \rbrace$.
By Lemma \ref{lemma estabilizadores finitos},  $\mathrm{Stab}_G(v)$ is a non-trivial $p$-subgroup of $\mathrm{SL}_2(k)$, since $G$ has not $p'$-torsion. 
Therefore, it follows from Lemma \ref{Lemma action of a p-group} that such a group fixes exactly one element $D_v$ in $\mathbb{P}^1(k)$. 

\begin{definition}\label{def of Q_x}
Let us denote by $Q_v=Q(v,D_v)$ the unique sector in $\X$ whose tip is $v$, and whose direction is $D_v$ introduced above.
\end{definition}

Since $\mathrm{Stab}_G(v)$ stabilizes both $v$ and $D_v$, it fixes $Q_v$ pointwise. 
The next result follows.

\begin{lemma}\label{Lemma stab in Q}
With the above notation, for each point $y \in Q_v$ we have $\mathrm{Stab}_G(v) \subseteq \mathrm{Stab}_G(y)$. In particular $Q_v \subseteq \Xinf$.
\end{lemma}

We are now in a position to prove Proposition \ref{lemma H_1}.

\begin{proof}[Proof of Proposition  \emph{\ref{lemma H_1}}] The result is consequence of the following claim: we assert that the map
\begin{equation}
   v \mapsto D_v,
\end{equation}
induces a bijection between the set of connected components of $\Xinf$ and $\mathbb{P}^1(k)$. 

Assume first that $v$ and $v'$ belong to the same connected component of $\Xinf$.
We must prove that $D_v=D_{v'}$.
Choosing an edge path from $v$ to $v'$, it suffices to treat the case in which $v$ and $v'$ are adjacent.
In this case, let $e$ be an edge in $\Xinf$ connecting $v$ with $v'$.
Then, $\mathrm{Stab}_G(e)$ is a non-trivial $p$-subgroup of $\mathrm{SL}_2(\OS)$, whence $\mathrm{Stab}_G(e)$ stabilizes one and only one vectorial chamber $D_e$ of $\mathbb{P}^{1}(k)$, by Lemma \ref{Lemma action of a p-group}.
Therefore, since $\mathrm{Stab}_G(e) \subseteq \mathrm{Stab}_G(v) \cap \mathrm{Stab}_G(v') $, we deduce that $\mathrm{Stab}_G(e)$ also stabilizes $D_v$ and $D_{v'}$. By uniqueness of $D_e$, we conclude $D_v=D_e=D_{v'}$.

Assume now $v,v' \in \Xinf$ satisfy $D_v=D_{v'}$. We claim that $v$ and $v'$ belong to the same connected component.
Indeed, $Q_v \cap Q_{v'}$ contains a common subsector by Lemma \ref{lemma equal visual limit}.
Since $Q_v$ and $Q_{v'}$ are both connected and have non-empty intersetion, the union $C=Q_v \cup Q_{v'}$ is connected. Since $v$ and $v'$ belong to $C\subseteq \Xinf$, the claim follows.

To prove surjectivity, fix $D\in\mathbb P^1(k)$ and choose
$g\in\mathrm{SL}_2(k)$ with $g\infty=D$. Write
$A=\mathcal O_{\mathcal S}$, $\Gamma=\mathrm{SL}_2(A)$, and let
$E_{12}$ be the elementary matrix with sole non-zero entry $1$ in
position $(1,2)$. Choose $0\ne d\in A$ such that
$d g E_{12}g^{-1}\in M_2(A)$. Then
\[
 V=\{I+xgE_{12}g^{-1}:x\in dA\}
\]
is an infinite subgroup of $\Gamma$, and every non-identity element of
$V$ is unipotent with unique fixed line $D$. Since
$[\Gamma:G\cap\Gamma]<\infty$, the subgroup $V\cap G$ has finite
index in $V$. Choose $1\ne u\in V\cap G$. The cyclic group generated
by $u$ has order $p$. In each tree $X_P$ it fixes a vertex: the convex
hull of a finite vertex orbit is a finite invariant subtree, and the
action has no edge inversions. Taking the product of these vertices
gives a vertex $v$ of $X$ fixed by $u$. Hence $v\in X_\infty$.
Since $u\in G_v$, it fixes $D_v$, and its unique fixed line is $D$.
Thus $D_v=D$.

Finally, $G_{av}=aG_va^{-1}$ gives $D_{av}=aD_v$ for $a\in G$.
Consequently $H_0(X_\infty;R)\cong R[\mathbb P^1(k)]$ as an
$R[G]$-modules. Under this identification the map to $H_0(X;R)=R$
is the augmentation. Since $X$ is contractible, the long exact
sequence of the pair identifies $H_1(X,X_\infty;R)$ with its kernel,
namely $\mathrm{St}_R$.
\end{proof}

\begin{corollary}\label{coro connected components as nerves}
Each connected component of $\Xinf$ has the form
$C_D:=\bigcup_{w\in I_D}Q_w$, where $I_D=\lbrace w\in\X_\infty\mid D_w=D \rbrace$.
\end{corollary}

\begin{proof}
We prove the result in two
steps. First that $C_D$ is a connected subset of $\Xinf$, and then that
it is maximal with this property.
We start with the first part. It follows from Lemma \ref{Lemma stab in Q} that $Q_w\subseteq \Xinf$.
Fix $w_0\in I_D$. By Lemma \ref{lemma equal visual limit}, for every $w\in I_D$ one has $
Q_w\cap Q_{w_0}\neq\varnothing$.
Since $Q_w$ and $Q_{w_0}$
are both connected (they are sectors, i.e., products of rays), and their
intersection is non-empty, the union $Q_w\cup Q_{w_0}$ is connected. 
Thus every $w \in C_D$ lies in the connected component of $w_0$, whence $C_D$ is connected.

Fix $v\in I_D$, and let $\mathcal C_v$ be the connected component of $v$ in $\Xinf$. Since $C_D$ is connected, $C_D\subseteq \mathcal C_v.$
Let $x\in \mathcal C_v$. 
Since $\Xinf$ is a simplicial complex
(in particular locally path-connected through its edges), there is an edge
path $v=v_0,\,v_1,\,\dots,\,v_n=x$
inside $\X$, with $v_i$ adjacent to $v_{i+1}$ for each $i$.
Since $D_{v_i}=D_{v_{i+1}}$ as shown in the proof of Lemma \ref{lemma H_1}, we have
$D_x=D_{v_n}=\cdots=D_{v_0}=D_v=D$. Hence $x\in I_D$, and since $x\in Q_x$ because $x$ is the tip of its own sector $Q_x$, we conclude $x\in C_D$, whence
$\mathcal C_v\subseteq C_D$ at the level of vertices. The statement for
higher-dimensional simplices of $\mathcal C_v$ follows by observing that such
a simplex has all of its vertices in $\mathcal C_v$, hence in $C_D$, and $C_D$
is a full subcomplex of $\Xinf$ by construction.
Finally, since any vertex of $\Xinf$ fixes an element $D\in \mathbb{P}^1(k)$, the result follows.
\end{proof}

We now prove that each connected component of $\Xinf$ is contractible. As in Corollary \ref{coro connected components as nerves}, we denote by $ I_D=\lbrace w\in\X_\infty\mid D_w=D \rbrace$
the set of vertices of $\Xinf$ whose associated direction is $D$, and $
C_D=\bigcup_{w\in I_D}Q_w.$

\begin{lemma}\label{lem:nerve-full}
Let $\mathcal N_D$ be the nerve of the cover $\{Q_w\}_{w\in I_D}$ of
$C_D$. Then $\mathcal N_D$ is the full simplex $\Delta^{I_D}$
on the vertex set $I_D$, i.e., every finite subset of $I_D$ spans a simplex of $\mathcal N_D$.
\end{lemma}

\begin{proof}
This is immediate from Corollary~\ref{cor:finite-int}. Indeed, by definition, a finite subset $F\subset I_D$ spans a simplex of $\mathcal N_D$ if and only if $\bigcap_{w\in F}Q_w\neq\varnothing$, and this holds for every finite family $F\subset I_D$, by Corollary~\ref{cor:finite-int}.
\end{proof}

\begin{lemma}\label{lem:full-simplex-contractible}
For any set $I$, the full simplex $\Delta^I$, equipped with
the weak topology as the colimit of its finite subsimplices, is contractible.
\end{lemma}

\begin{proof}
Fix $w_0\in I$ and define $H:\Delta^I\times[0,1]\to\Delta^I$ via barycentric coordinates as follows: if $x=\sum_{w\in I}t_w\, w$ is a finite sum with $t_w\ge0$ and $\sum
t_w=1$), set:
\[
H(x,s) := (1-s)\,x + s\, w_0 = \sum_{w\in I} (1-s)t_w\, w \;+\; s\, w_0.
\]
This is precisely the cone structure with apex $w_0$ on $\Delta^I$.
Each point $x$ has finite support $F=\{w: t_w>0\}$, so $H(x,s)$ has support contained in $F\cup\{w_0\}$, which is finite, hence $H(x,s)\in\Delta^{F\cup\{w_0\}}\subset
\Delta^I$ for all $s$. The map $H$ is continuous with respect to the weak
(colimit) topology, since its restriction to each finite subsimplex
$\Delta^F\times[0,1]$ is continuous (it is affine in barycentric coordinates),
and the weak topology is exactly the final topology with respect to these
inclusions. Moreover $H(x,0)=x$ and $H(x,1)=w_0$ for all $x$. Hence $H$ is a
homotopy between $\mathrm{id}_{\Delta^I}$ and the constant map $w_0$, i.e.,
$\Delta^I$ is contractible. $\qed$
\end{proof}

\begin{proposition}\label{prop:CD-simply-connected}
For every $D\in\mathbb P^1(k)$, the connected component $C_D \subseteq \Xinf $, as defined in Corollary \ref{coro connected components as nerves}, is contractible. 
In particular, every connected component of $\Xinf$ is simply connected.
\end{proposition}

\begin{proof}
The cover $\{Q_w\}_{w\in I_D}$ of $C_D$ satisfies the hypotheses of Lemma \ref{lem:nerve-full}. Each $Q_w$ is a contractible subcomplex of
$C_D$ since it is a product of rays.
Moreover, every non-empty finite intersection $Q_{w_1}\cap\cdots\cap Q_{w_n}$ is again sector chamber by Corollary~\ref{cor:finite-int}. Hence this is contractible. By
the nerve lemma, it follows that $
C_D \cong |\mathcal N_D|.$
By Lemma~\ref{lem:nerve-full}, $\mathcal N_D=\Delta^{I_D}$, and by
Lemma~\ref{lem:full-simplex-contractible}, $\Delta^{I_D}$ is contractible.
Therefore $C_D$ is contractible.
Since every connected component of $\Xinf$ is of the form
$C_D$ for some $D\in\mathbb P^1(k)$ by Corollary \ref{coro connected components as nerves}, the result follows.
\end{proof}

\section{A resolutions of \texorpdfstring{$\mathrm{St}$ and its flat dimension}{St}}\label{section res and dim st}

We begin this section by showing that the Steinberg module $\mathrm{St}$ introduced in \S \ref{section Stainbert} is not flat whenever $s \geq 2$. The case of $s=1$ has been covered by Serre in \cite{SerreTrees}, showing that $\mathrm{St}$ is projective, so in particular flat.

Let $s \geq 2$. We will prove that $\mathrm{St}$ is not flat over $R[G]$. Indeed, by Dirichlet's unit theorem we have $\mathrm{rk}_{\mathbb{Z}}\left(\mathcal{O}_S^\times\right) = s - 1 > 0$.
Let $a \in \mathcal{O}_S^\times \smallsetminus \mathcal{O}_{S,\mathrm{tor}}^\times$, and let $\tau \in \mathrm{SL}_2(\mathcal{O}_S)$ be given by
\[\tau = \begin{pmatrix}
        a & 0 \\ 0 & a^{-1}
\end{pmatrix}.\]
Since $G$ has finite index in $\mathrm{SL}_2(\mathcal{O}_S)$, there exists a positive integer $n \ge 1$ such that $\tau^n \in G$. By replacing $\tau$ with $\tau^n$ if necessary, we may assume without loss of generality that $\tau \in G$. 
Let $T$ be the subgroup of $G$ generated by $\tau$.
The element $\tau$ acts on $\mathbb{P}^1(k)$ via $\tau \cdot z = a^2 z$. This action fixes exactly two points, $P = 0$ and $Q = \infty$, and acts freely on $\mathbb{P}^1(k) \setminus \{0, \infty\}$. Consequently, we have
\[R\left[\mathbb{P}^1(k)\right]^T = R[P] \oplus R[Q] \cong R^2.\]
Applying the $T$-invariants functor to the short exact sequence
\[0 \to \mathrm{St} \to R\left[\mathbb{P}^1(k)\right] \to R \to 0,\]
yields the exact sequence
\[0 \to \mathrm{St}^T \to R[P] \oplus R[Q] \xrightarrow{\mathrm{aug}} R.\]
Since $\mathrm{aug}(P) = \mathrm{aug}(Q) = 1$, the augmentation map $\mathrm{aug}$ is surjective, and its kernel is $\mathrm{St}^T$. Thus, we obtain $\mathrm{St}^T = \langle [P] - [Q] \rangle$.

On the other hand, consider the short exact sequence 
\[0 \to R[T] \xrightarrow{\tau-1} R[T] \to R \to 0.\]
Tensoring with $\mathrm{St}_R$, we obtain the long exact sequence
\[\cdots \to 0 \to \operatorname{Tor}_1^{R[T]}(R, \mathrm{St}_R) \to \mathrm{St}_R \xrightarrow{\tau-1} \mathrm{St}_R \to \cdots.\]
This yields $\operatorname{Tor}_1^{R[T]}(R, \mathrm{St}_R) \cong \mathrm{St}_R^T \neq 0$. 
In particular, $\operatorname{Tor}_1^{R[T]}(R, \mathrm{St}_R) \neq 0$.
\begin{comment}
On the other hand, note that $\mathrm{Tor}_1^{R[T]}(R, \mathrm{St}) \cong \mathrm{St}^T \neq 0$.
\green{Explicar este isomorfismo.}
\end{comment}
Therefore, the restriction $\mathrm{Res}_T(\mathrm{St})$ is not $R[T]$-flat. Since $R[G]$ is a free $R[T]$-module (and flat base change preserves flatness), it follows that $\mathrm{St}$ cannot be a flat $R[G]$-module.\\

In what follows we describe a (free, flat and projective) resolution of $\mathrm{St}$ in terms of the $i$-dimensional relative chains $L_i$ of $X$ modulo $X_\infty$ with coefficients in $R$. Indeed, we have the complex of relative chains:

\begin{equation}\label{sucesion exacta de cadenas relativas}
    0\to L_s \to L_{s-1}\to \cdots \to L_1\to L_0 \to 0.
\end{equation}

Let $B_i$ and $Z_i$ be the relative boundaries and relative cycles, respectively. Since $H_{i}(X\, \mathrm{mod }\, X_\infty;R)$ is trivial for all $i>1$, the complex in \eqref{sucesion exacta de cadenas relativas} induces the following exact sequences:

\begin{equation}\label{sucesion exacta para B1}
    0\to L_s \to L_{s-1}\to \cdots \to L_2\to B_1 \to 0.
\end{equation}
and
\begin{equation}\label{sucesion exacta para Z1}
    0\to Z_1 \to L_1\to  L_0 \to 0.
\end{equation}

\begin{lemma}\label{lem splitting of L}
The complex $L_\bullet$ defined in \eqref{sucesion exacta de cadenas relativas} splits as a
direct sum of complexes of $R[G]$-modules
$L_\bullet\ =\ P_\bullet[1]\ \oplus\ C_\bullet,$
where:
\begin{itemize}
  \item[(i)] $C_\bullet$ is the contractible complex $0\to L_0\xrightarrow{\ \mathrm{id}\ }L_0\to 0$
        concentrated in degrees $1$ and $0$;
  \item[(ii)] $P_\bullet$ is the complex
        $0\to P_{s-1}\to \cdots \to P_1 \to P_0\to 0,$ where thr $P_j$ are finitely generated projective $R[G]$-modules given by $P_j:=L_{j+1}$ for $j\geq 1$ and $ P_0:=Z_1$.
        The complex $P_\bullet[1]$ denotes
        $P_\bullet$ shifted in $1$, so that $P_j$ sits in degree $j+1$.
\end{itemize}
\end{lemma}
 
\begin{proof}
Since $\X$ is connected and $\Xinf\neq\emptyset$, the map
$H_0(\Xinf ;R)\to H_0(\X ;R)=R$ is surjective, so the long exact sequence
\eqref{eq top rel homology} gives $H_0(\X\bmod\Xinf;R)= \lbrace 0 \rbrace$. 
Therefore, the differential
$d_1\colon L_1\to L_0$ is surjective. As $L_0$ is free, $d_1$ admits a $R[G]$-linear
section $\varsigma\colon L_0\to L_1$, whence
$L_1=Z_1\oplus \varsigma(L_0)$ and $Z_1=\ker(d_1)$.
In particular $Z_1$ is a direct summand of the finitely generated free module $L_1$,
hence finitely generated projective.
Put $C_\bullet:=\big(0\to\varsigma(L_0)\xrightarrow{\,d_1\,}L_0\to 0\big)$, which is a
subcomplex of $L_\bullet$ because $d_1$ restricted to $\varsigma(L_0)$ is an isomorphism
onto $L_0$. On the other hand, $d_1\circ d_2=0$ forces
$d_2(L_2)\subseteq Z_1$ and $d_1|_{Z_1}=0$, so the graded submodule with $L_j$ in
degree $j\geq 2$ and $Z_1$ in degree $1$ is also a subcomplex, namely $P_\bullet[1]$.
The two subcomplexes are complementary in every degree, which proves the claim.
\end{proof}

We now give an upper bound for the projective and flat dimensions of the Steinberg module $\mathrm{St}$.

\begin{proposition}\label{prop P is a resolution}
The complex $P_\bullet$ of Lemma \ref{lem splitting of L}, augmented by the canonical
surjection $P_0=Z_1\twoheadrightarrow Z_1/B_1$, is a resolution of $\St$ by finitely
generated projective and flat $R[G]$-modules of length $s-1$.
In particular both
$\mathrm{pd}_{R[G]}(\St)$, $\mathrm{fd}_{R[G]}(\St)  \leq\ s-1 $.
\end{proposition}
 
\begin{proof}
By Proposition \ref{lemma H_1} together with
Proposition \ref{prop:CD-simply-connected} we have
$H_j(\X\bmod\Xinf; R)=\St$ for $j=1$ and $H_j(\X\bmod\Xinf;R )= \lbrace 0 \rbrace$ for $j\neq 1$. Since
$H_j(P_\bullet)=H_{j+1}(L_\bullet)$ for $j\geq 1$ and
$H_0(P_\bullet)=Z_1/B_1=H_1(L_\bullet)=\St$, the augmented complex
$P_\bullet\to\St\to 0$ is exact. Each $P_j$ is finitely generated, projective and flat by Lemma \ref{lem splitting of L}, and $P_j=0$ for $j>s-1$.
\end{proof}

In the sequel we find a lower bound for the flat dimension of $\mathrm{St}$. To do so, we use the next lemma.

\begin{lemma}\label{lemma torus rank s}
The group $G$ contains a subgroup isomorphic to $\mathbb{Z}^{s}$.
\end{lemma}

\begin{proof}
We first choose a separable quadratic extension $\ell/k$ in which every $P \in \mathcal{S}$
splits. Fix a closed point $Q \notin \mathcal{S}$. If $p=\mathrm{char}(\mathbb{F}) \neq 2$, the strong approximation theorem
gives $f \in k$ with $\nu_P(f-1) \geq 1$ for all $P \in \mathcal{S}$ and $\nu_Q(f) = 1$ for $Q \notin \mathcal{S}$.
Then $f$ is a square in each $k_P$ by Hensel's lemma, but not in $k$, and we set
$\ell = k(\sqrt{f})$.
If $p = 2$, choose instead $f \in k$ with $\nu_P(f) \geq 1$ for all
$P \in \mathcal{S}$ and $\nu_Q(f) = -1$ for $Q \notin \mathcal{S}$. Then $x^2 - x = f$ has a root in each $k_P$ by
Hensel's lemma, but none in $k$, and we let $\ell$ be the corresponding
Artin-Schreier extension.

Let $\mathcal{O}_\ell=\mathcal{O}_{\mathcal{S},\ell}$ be the integral closure of $\mathcal{O}_\mathcal{S}$ in $\ell$. Since every
$P \in \mathcal{S}$ splits in $\ell$, there are exactly $2s$ places of $\ell$ above $\mathcal{S}$, so
$\mathcal{O}_\ell^{*}$ has $\mathbb{Z}$-rank $2s-1$ by Dirichlet's unit theorem. 
The image of the
norm $N_{\ell/k} \colon \mathcal{O}_\ell^{*} \to \mathcal{O}_{\mathcal{S}}^{*}$ contains
$(\mathcal{O}_\mathcal{S}^{*})^2$, so it has $\mathbb{Z}$-rank $s-1$, that is maximal by Dirichlet's unit theorem.
Hence $\mathcal{T} := \ker(N_{\ell/k}) \cap \mathcal{O}_\ell^{*}$ has rank $s$.

The group $\mathcal{T}$ acts by multiplication on $\mathcal{O}_\ell$, which is a projective $\mathcal{O}_\mathcal{S}$-module of rank $2$ (over $\OS$), with determinant $N_{\ell/k} = 1$.
Choosing a $k$-basis of $\ell$, this identifies $\mathcal{T}$ with a subgroup of
$\mathrm{SL}(\mathcal{O}_\ell) \subseteq \mathrm{SL}_2(k)$. Since the stabilizers of two
$\OS$-lattices of $k^2$ are commensurable, $\mathrm{SL}(\mathcal{O}_\ell)$ is
commensurable with $\mathrm{SL}_2(\OS)$, hence with $G$. Therefore
$\mathcal{T} \cap G$ has finite index in $\mathcal{T}$. It is a finitely generated
abelian group of rank $s$, so it contains a subgroup isomorphic to $\mathbb{Z}^s$.
\end{proof}

\begin{theorem}\label{theorem pd}
The projective and flat dimension of $\mathrm{St}$ over $R[G]$ is $s-1.$
\end{theorem}

\begin{proof}
Let $T' \subseteq G$ be the subgroup isomorphic to $\mathbb{Z}^{s}$, given by
Lemma~\ref{lemma torus rank s}. Each $t \in T' \smallsetminus \{1\}$ is an element of
$\ell \smallsetminus k$, since, by construction, the only elements of $\mathcal{T} \cap k$ are $\pm 1$ and $T'$ is
torsion-free.
Hence the eigenvalues $t, \bar{t}$ of $t$ acting on $\ell \cong k^2$ do not
lie in $k$, so $t$ fixes no point of $\mathbb{P}^1(k)$. Thus $T'$ acts freely on
$\mathbb{P}^1(k)$, and $R[\mathbb{P}^1(k)]$ is a free $R[T']$-module. Applying
$H_*(T', -)$ to the short exact sequence
\[0\xrightarrow[]{} \mathrm{St_R}\xrightarrow[]{}R\left[ \mathbb{P}^1(k)\right]\xrightarrow[]{}R\xrightarrow[]{}0,\]
and using that $H_s\big(T', R[\mathbb{P}^1(k)]\big) = 0$ since $s \geq 1$, we
obtain an injection $R \cong H_s(T', R) \hookrightarrow H_{s-1}(T', \mathrm{St})$, by
Koszul's lemma.
On the other hand, by Shapiro's lemma,
$$ H_{s-1}(T', \mathrm{St}) = \mathrm{Tor}_{s-1}^{R[T']}(R,\mathrm{St})\simeq \mathrm{Tor}_{s-1}^{R[G]}(R\otimes_{R[T']} R[G],\mathrm{St}).$$
Combining these isomorphisms yields
$\mathrm{Tor}_{s-1}^{R[G]}(R\otimes_{R[T']} R[G],\mathrm{St})\neq \{0\}.$
In particular, the flat dimension of $\mathrm{St}$ is at least $s-1$.
Since flat dimension is bounded above by projective dimension, we conclude that
$\mathrm{pd}_{R[G]}(\mathrm{St}) \ge s - 1$.
Hence the result follows from Prop.~\ref{prop P is a resolution}.
\end{proof}

\begin{remark}
We can alternatively prove the preceding theorem by using the diagonal $\mathrm{T}(k)\subset \mathrm{SL}_2(k)$. However, we use the group $\mathbb{Z}^{s}$ by uniformity with the subsequent results.
Define $\mathrm{T}(\mathcal{O}_S)=\mathrm{T}(k)\cap \mathrm{SL}_2(\mathcal{O}_s)$ and set $A=T(\mathcal{O}_S)\cap G.$
Since $\mathrm{T}(\mathcal{O}_s)\simeq \mathcal{O}_S^\times$ and $G$ is $p'$-torsion free, it follows that $A\simeq \mathbb{Z}^{s-1}$. The subgroup $A$ acts on $\mathbb{P}^1(k)$ fixing $0$ and $\infty$. Applying the $A$-invariants functor to the short exact sequence
\[0\xrightarrow[]{} \mathrm{St_R}\xrightarrow[]{}R\left[ \mathbb{P}^1(k)\right]\xrightarrow[]{}R\xrightarrow[]{}0,\] we obtain $\mathrm{St}_R^{A}\simeq R.$
By Koszul's lemma, we have $\mathrm{Tor}_{s-1}^{R[A]}(R,\mathrm{St})\simeq \mathrm{St}^{A}_R\simeq R.$ 
On the other hand, by Shapiro's lemma,
\[\mathrm{Tor}_{s-1}^{R[A]}(R,\mathrm{St})\simeq \mathrm{Tor}_{s-1}^{R[G]}(R\otimes_{R[A]} R[G],\mathrm{St}).\]
Combining these isomorphisms yields $\mathrm{Tor}_{s-1}^{R[G]}(R\otimes_{R[A]} R[G],\mathrm{St})\neq \{0\}.$ 
In particular, the flat dimension of $\mathrm{St}$ is at least $s-1$. 
Therefore, we conclude as in the preceding proof.
\end{remark}

\section{Relative Homology and the Euler-Poincar\'e characteristic}\label{section Rel. Homology}

We are now in position to compute the relative homological dimension of $G$ modulo $\{G_\sigma\}_{\sigma \in \Sigma}$ introduced in Def. \ref{def relative homology}.
In particular, here we prove Ths. \ref{main theorem intro pd}, \ref{main theorem intro chi}, \ref{main theorem intro euler}.
We start by proving the first one, with is a direct consequence of Th. \ref{theorem pd}.

\begin{proof}[Proof of Th. \ref{main theorem intro pd}]
It follows from Theorem \ref{theorem pd} that $\mathrm{rhd}_{R}(G,G_{\sigma})=\mathrm{fd}_{R[G]}(\mathrm{St})+1=s$. 
\end{proof}

The next two result are consequences of Th. \ref{main theorem intro pd}.

\begin{corollary}
For every $r>s$ and every $R[G]$-module $M$, we have $H_r(G \bmod G_{\sigma},M)=\lbrace 0\rbrace$. \qed
\end{corollary}

The $\mathrm{Tor}$ exact sequence, together with Shapiro's lemma, give the exact sequence:
$$\cdots \to H_{i+1}( G \text{ mod } G_{\sigma}, M) \to 
\bigoplus_{\sigma \in \Sigma_0} H_i(G_{\sigma},M) \to
H_i (G,M) \to \cdots .$$
Therefore, the next result immediate.

\begin{corollary}\label{cor direct sum}
For each $r >s$ and every $R[G]$-module $M$, we have $H_{r}(G,M) \cong \bigoplus_{\sigma \in \Sigma_0} H_r(G_{\sigma},M)$. For $r=s$, we have the exact sequence
$0 \to \bigoplus_{\sigma \in \Sigma_0} H_s(G_{\sigma},M) \to H_s (G,M) \to H_s(G \bmod G_{\sigma},M).$ \qed
\end{corollary}

Keeping the notation of
\eqref{sucesion exacta de cadenas relativas},
recall that $G$ acts freely on each $S_i$, so that a choice of orbit representatives yields an isomorphism of $R[G]$-modules $L_i\simeq R[G]^{l_i}$, which is finitely generated by Prop. \ref{prop li finite}.

\begin{lemma}\label{thm comparison}
Let $M$ be an $R[G]$-module. Then, for every $n\in\mathbb{Z}_{\geq 0}$, there is a natural
isomorphism of abelian groups
$$H_n\big(L_\bullet\otimes_{R[G]}M\big)\ \cong\ H_n\big(G\bmod G_\sigma,M\big)
\ =\ \mathrm{Tor}^{R[G]}_{n-1}(\St,M),$$
with the convention $\mathrm{Tor}^{R[G]}_{-1}=0$. In other words, the relative homology
of $G$ modulo $\lbrace G_\sigma\rbrace_{\sigma\in\Sigma}$ is the $G$-equivariant
homology of the pair $(\X,\Xinf)$ with coefficients in $M$.
\end{lemma}
 
\begin{proof}
Tensoring the decomposition of Lemma \ref{lem splitting of L} with $M$ over $R[G]$ gives
$$L_\bullet\otimes_{R[G]}M\ =\ \big(P_\bullet\otimes_{R[G]}M\big)[1]\ \oplus\ \big(C_\bullet\otimes_{R[G]}M\big).$$
The complex $C_\bullet\otimes_{R[G]}M$ is again of the form
$0\to M'\xrightarrow{\ \cong\ }M'\to 0$, hence acyclic, so it does not contribute to
homology. Therefore
$H_n\big(L_\bullet\otimes_{R[G]}M\big)\ \cong\ H_{n-1}\big(P_\bullet\otimes_{R[G]}M\big).$
By Proposition \ref{prop P is a resolution}, $P_\bullet$ is a projective resolution of
$\St$, so the right-hand side is $\mathrm{Tor}^{R[G]}_{n-1}(\St,M)$ by definition of the
derived functor \cite[Ch. III]{Brown-Cohomology}. Finally, by
Lemma \ref{lemma Deleta=St} the module $\Delta=\ker(\varepsilon)$ of
Definition \ref{def relative homology} is exactly $\St$, whence
$\mathrm{Tor}^{R[G]}_{n-1}(\St,M)=H_n(G\bmod G_\sigma,M)$.
\end{proof}

We now focus on describing the relative homology groups $H_r(G \text{ mod } G_\sigma, M)$, for $r=1,...,s$ in terms of the Euler-Poincar\'e characteristic $\chi(G)$ of $G$.
We start by proving Th. \ref{main theorem intro chi}.

\begin{proof}[Proof of Th. \ref{main theorem intro chi}]
We decompose the sum introduced in Eq. \eqref{eq euler-poincare char} in two parts. One relative to $\X_{\infty}$ and other relative to $\X$ mod $\Xinf$.
When $F$ is a simplex outside $\Xinf$ we have $\mathrm{Stab}_G(F)= \lbrace e \rbrace$.
Then, the contribution of $\X \text{ mod } \Xinf$ to $\chi(G)$ is exactly
$$\chi(G)|_{\X \text{ mod } \Xinf } = \sum_{i=0}^{s}(-1)^{i} \#(G \backslash S_i) = \sum_{i=0}^{s} (-1)^{i} l_i .$$
Therefore, in order to conclude the result, we just need to show that the contribution of $\X_{\infty}$ to $\chi(G)$ is zero.
In other words, we have to prove that
\begin{equation}\label{eq chi in X}
    \chi(G)|_{\Xinf}=  \sum_{F \in G\backslash \Xinf}  \frac{(-1)^{\mathrm{dim}(F)}}{ \#\mathrm{Stab}_G(F)} =0.
\end{equation}
%\begin{equation}
  %  \sum_{F \in G\backslash X_{\infty}} \frac{1}{ \mathrm{Card}(\mathrm{Stab}_G(F))} + \cdots + \sum_{F \in G\backslash X_{\infty}} \frac{1}{ \mathrm{Card}(\mathrm{Stab}_G(F))} + \cdots \sum_{F \in G\backslash X_{\infty}} \frac{1}{ \mathrm{Card}(\mathrm{Stab}_G(F))}  =0. 
%\end{equation}
Let $v$ be a vertex in $\Xinf$. Let $Q_v=Q(v,D_v)$ be the sector introduced in Definition \ref{def of Q_x}.
Let $C=C_v$ be the unique alcove of $Q_v$ that contains $v$.
Since $\mathrm{Stab}_G(v)$ stabilizes $Q_v$ and $G \subseteq \mathrm{SL}_2(k)$ does not change the type of vertices, we get
$$\mathrm{Stab}_G(v) \subseteq \mathrm{Stab}_G(C_v) = \mathrm{Fix}_G(C_v). $$
Thus, we deduce $\mathrm{Stab}_G(v)=\mathrm{Stab}_G(C_v)$, because always $\mathrm{Stab}_G(C_v) \subseteq \mathrm{Stab}_G(v)$.
In particular, for each facet $F$ satisfying $\lbrace v \rbrace \subseteq F \subseteq C_v$ we have
$$\mathrm{Stab}_G(v)=\mathrm{Stab}_G(F)=\mathrm{Stab}_G(C_v). $$
Thus, the contribution of $v$ and all the facets of $C_v$ containing $v$ to $\chi(G)$ cancel out. 
Therefore:
$$ 
\chi(G)|_{\Xinf}= \sum_{\bar{v} \in G \backslash \Xinf} \sum_{\bar{v} \in \bar{F}} \frac{(-1)^{\mathrm{dim}(F)}}{ \#\mathrm{Stab}_G(F)} =0.
$$
\end{proof}

Let $D$ be a principal ideal domain, and let $F$ be its fraction field. For all $D$-modules $A$, we write $\rk_{D}(A):=\dim_{F}(A\otimes_{D}F)
\in \mathbb{Z}_{\geq0}\cup\lbrace\infty\rbrace$. Since $-\otimes_{D}F$ is exact,
$\rk_{D}$ is additive on short exact sequences of $D$-modules and,
consequently, the Euler characteristic of a bounded complex of finite rank $D$-modules equals the Euler characteristic of its homology. 
The next result is exactly Th. \ref{main theorem intro euler}.

\begin{theorem}\label{thm euler relative}
Let $M$ be a $D[G]$-module of finite rank
$m:=\rk_{D}(M)<\infty$. Then:
\begin{itemize}
  \item[(a)] $H_i(G\bmod G_\sigma,M)$ has finite rank for every $i$, and vanishes for
        $i>s$ and for $i=0$, and
  \item[(b)] one has the identity
        $$\sum_{i=1}^{s}(-1)^{i} \rk_{D} H_i\big(G\bmod G_\sigma,M\big) =\chi(G)\cdot m .$$
\end{itemize}
\end{theorem}

\begin{proof}[Proof of Th. \ref{main theorem intro euler}]
We have $L_i\cong D[G]^{\,l_i}$, whence an isomorphism of abelian
groups $L_i\otimes_{D[G]}M \cong M^{l_i}$ so that $\rk_{D}\big(L_i\otimes_{D[G]}M\big)=l_i m<\infty .$
Thus $L_\bullet\otimes_{D[G]}M$ is a complex of abelian groups of finite rank,
concentrated in degrees $0,\dots,s$. Its homology groups are subquotients of its terms,
hence of finite rank, and vanish outside $\lbrace 0,\dots,s\rbrace$. By
Theorem \ref{thm comparison} these homology groups are the
$H_i(G\bmod G_\sigma,M)$. Moreover $H_0(G\bmod G_\sigma,M)=0$ by convention. This proves statement (a).
 
We now prove statement (b). Note that the Euler characteristic of the bounded complex $L_\bullet\otimes_{D[G]} M$ may be
computed either from its terms or from its homology.
More concretely, we have
$$\sum_{n=0}^{s}(-1)^{n} \rk_{D}\big(L_n\otimes_{D[G]}M\big) = \sum_{n=0}^{s}(-1)^{n} \rk_{D} H_n\big(L_\bullet\otimes_{D[G]}M\big).$$
The left-hand side equals $\big(\sum_{n=0}^{s}(-1)^{n}l_n\big)m=\chi(G)\,m$ by
Th. \ref{main theorem intro chi}, and by Th. \ref{thm comparison} the
right-hand side equals $\sum_{n=1}^{s}(-1)^{n}\rk_{D}H_n(G\bmod G_\sigma,M)$, whence the result follows.
\end{proof}
 
\begin{corollary}\label{cor euler trivial coeff}
Taking $M=D$ with the trivial $G$-action, we have
\[\sum_{i=1}^{s}(-1)^{i} \rk_{D} H_i\big(G\bmod G_\sigma,D\big) = \chi(G).\] \qed
\end{corollary}

\begin{example}
For $s=1$ one recovers the classical result of \cite[Ch. II]{SerreTrees}. There
$\mathrm{fd}_{\ZG}(\St)=0$, so $H_i(G\bmod G_\sigma,M)=0$ for $i\neq 1$ and
Theorem \ref{thm euler relative} reduces to
$\rk_{\mathbb{Z}}\big(\St\otimes_{\ZG}M\big)=-\chi(G) m$.
\end{example}

\begin{remark}\label{rem only alternating}
The Euler characteristic sees the alternating sum and nothing else, so
Theorem \ref{thm euler relative} does not by itself determine the individual ranks. It
does so as soon as one knows a concentration statement. If
$H_i(G\bmod G_\sigma,M)=0$ for all $i\neq i_0$, then
$\rk_{D}\,H_{i_0}\big(G\bmod G_\sigma,M\big)=(-1)^{i_0}\chi(G)\cdot m .$
Note that the module used in Theorem \ref{theorem pd} to bound
$\mathrm{fd}_{D[G]}(\St)$ from below is the induced module
$D \otimes_{D[T']}D[G]$, which has infinite rank over $D$ and hence lies
outside the range of Theorem \ref{thm euler relative}. It is therefore consistent to ask
whether $\mathrm{Tor}^{D[G]}_j(\St,M)=0$ for $j\geq 1$ whenever $\rk_{D}(M)<\infty$,
equivalently whether the relative homology is concentrated in degree $1$ for coefficients
of finite rank.
\end{remark}

Let $h := \#\Sigma_0 = \#G\backslash\mathbb{P}^1(k)$ be the number of cusps of $G$. This number is finite. Indeed, when $G \subseteq \mathrm{SL}_2(\mathcal{O}_S)$ we have $h \leq \#\mathrm{Pic}(\mathcal{O}_S)\cdot[\mathrm{SL}_2(\mathcal{O}_S):G]$, and in general finiteness follows by applying this bound to the finite-index subgroup $G \cap \mathrm{SL}_2(\mathcal{O}_S)$, which has at least as many cusps as $G$.
Finally, $\#\mathrm{Pic}(\mathcal{O}_S)$ is finite by Weil's theorem, since $\mathbb{F}$ is finite.
%The previous number is finite by Weil's theorem, since $\mathbb{F}$ is finite.

\begin{corollary}\label{cor lower bound H1}
One has $\rk_{D} H_1\big(G\bmod G_\sigma,D\big)\ \geq\ h-1$.
\end{corollary}
 
\begin{proof}
Take $M=D$ in the long exact sequence
$$\cdots\to H_1\big(G\bmod G_\sigma,D\big)\to
\bigoplus_{\sigma\in\Sigma_0}H_0(G_\sigma,D)\to H_0(G,D)\to
H_0\big(G\bmod G_\sigma,D\big)=0 .$$
Since $H_0(G_\sigma,D)=D$ and $H_0(G,D)=D$, the second map is the
augmentation $D^{h}\to D$, whose kernel has rank $h-1$. By exactness, this kernel
is a quotient of $H_1(G\bmod G_\sigma,D)$.
\end{proof}

\section{Relative homology modulo maximal unipotent subgroups}\label{section Rel. Homology for U}

For each $\sigma \in \Sigma_0$, let $U_\sigma$ be the subgroup of unipotent
elements of $G_\sigma$.
As in Def. \ref{def relative homology}, let $\Delta_{\mathfrak{U}} $ be the kernel of the augmentation homomorphism $\bigoplus_{\sigma\in\Sigma_0} R[G/U_{\sigma}]\to R$, so that $H_i(G \bmod U_{\sigma}, M)= \mathrm{Tor}_{i-1}^{R[G]}(\Delta_{\mathfrak{U}},M)$, for all $i\geq 1$.
In this section, we compare the relative homology groups $H_*(G \bmod G_\sigma, M)$
with $H_*(G \bmod U_\sigma, M)$.
In particular, here we prove Th. \ref{main theo intro pd of delta} which states that flat dimension of $\Delta_{\mathfrak{U}}$ is exactly $s-1$.
After that, we describe the homology groups $H_r(G \bmod U_\sigma, M)$, for $r\leq s$, in terms of the Euler-Poincar\'e characteristic of $G$, by proving Th. \ref{main thm euler-unipotent for U}.
To do so, we need some preparation, which divides in next two $4$ results, the first of one relates the group structure of $G_{\sigma}$ and $U_{\sigma}$.

\begin{lemma}\label{lemma structure Gs}
The group $U_{\sigma}$ is a maximal unipotent subgroup of $G$, and the family
$\mathfrak{U}=\lbrace U_{\sigma} \rbrace_{\sigma \in \Sigma_0}$ is a representative system
of the $G$-conjugacy classes of maximal unipotent subgroups in $G$.
If the torsion of $G$ is $p$-primary, we have an exact sequence of the form $1 \to U_{\sigma} \to G_\sigma \to \mathbb{Z}^{s-1} \to 1$,
for all $\sigma \in \Sigma_0$.
\end{lemma}

\begin{proof}
The first two statements follow from \cite[Th.~3.1]{BravoLoisel2}.
Thus, we focus on proving the third one.
Let $\gamma_\sigma \in \mathrm{SL}_2(k)$ with $\gamma_\sigma \cdot \sigma = \infty$. By
\cite[Lemma~8.1]{BravoLoisel2}, we have $G_\sigma / U_\sigma \cong T(G_\sigma)$, where
$T(G_\sigma)$ is the image of $G_\sigma$ in
$\gamma_\sigma^{-1}\mathrm{T}(k)\gamma_\sigma \cong k^\times$. Since
$\gamma_\sigma G\gamma_\sigma^{-1}$ is commensurable with
$\mathrm{SL}_2(\mathcal{O}_S)$, the group
$\gamma_\sigma G_\sigma \gamma_\sigma^{-1} = \gamma_\sigma G\gamma_\sigma^{-1} \cap \mathrm{B}(k)$
is commensurable with $\mathrm{B}(\mathcal{O}_S)$, where $\mathrm{B}$ is the group of
upper triangular matrices in $\mathrm{SL}_2$.
Hence $T(G_\sigma)$ is commensurable with the image of $\mathrm{B}(\mathcal{O}_S)$,
namely $\mathcal{O}_S^\times$, so it has rank $s-1$.
It remains to see that $T(G_\sigma)$ is torsion-free. The image in $k^\times$ of a
torsion element $t \in T(G_\sigma)$ is a root of unity, hence lies in $\mathbb{F}^\times$
since $\mathbb{F}$ is algebraically closed in $k$. In particular, $t$ has
order $n$ prime to $p$. If $g = tu \in G_\sigma$ is a lift of $t$, then $g^n$ is unipotent, so $g$ has finite order, which is a power of $p$ by hypothesis.
Since $n$ divides this order, we have $t = 1$. Thus
$T(G_\sigma) \cong \mathbb{Z}^{s-1}$ is free, and we have the exact sequence $1 \to U_\sigma \to G_\sigma \to \mathbb{Z}^{s-1} \to 1$.
\end{proof}

\begin{lemma}\label{lemma U-acyclic}
Assume that $p \in R^{*}$. Then, for every $\sigma \in \Sigma_0$ and every
$R[U_\sigma]$-module $M$, we have $H_j(U_\sigma, M) = 0$ for all $j \geq 1$.
\end{lemma}

\begin{proof}
By \cite[Th.~3.1]{BravoLoisel2}, the group $U_\sigma$ is contained in $\gamma_\sigma^{-1}\mathrm{U}^+(k)\gamma_\sigma \cong (k, +)$, where
$\gamma_\sigma \in \mathrm{SL}_2(k)$ as in the proof of Lemma~\ref{lemma structure Gs}. Thus $U_\sigma$
is an elementary abelian $p$-group, hence the directed union of its finite subgroups.
For each finite subgroup $F \subseteq U_\sigma$ and each $j \geq 1$, the group
$H_j(F, M)$ is annihilated by $\# F$, which is invertible in $R$, so
$H_j(F, M) = 0$. Since homology commutes with directed colimits, the result
follows.
\end{proof}

For a group $H$ and an $R[H]$-module $M$, we write
$M_H = H_0(H, M)= M / \langle hm - m \vert h \in H, m \in M \rangle$ for the module
of coinvariants. 

\begin{corollary}\label{coro hochschild serre}
Assume that $p \in R^{*}$ and that the torsion of $G$ is $p$-primary. Then, for
every $R[G_\sigma]$-module $M$, we have
$H_*(G_{\sigma},M) \cong H_*(\mathbb{Z}^{s-1},M_{U_{\sigma}})$.
\end{corollary}

\begin{proof}
It follows from Lemma \ref{lemma structure Gs} that $1 \to U_{\sigma} \to G_\sigma \to \mathbb{Z}^{s-1}\to 1$ is exact.
Consider the Lyndon--Hochschild--Serre spectral sequence
$$
E^2_{i,j} = H_i\big(\mathbb{Z}^{s-1}, H_j(U_\sigma, M)\big) \Longrightarrow H_{i+j}(G_\sigma, M),$$
where $\mathbb{Z}^{s-1} \cong G_\sigma / U_\sigma$ acts on $H_j(U_\sigma, M)$, and in
particular on $H_0(U_\sigma, M) = M_{U_\sigma}$, via the induced action. By
Lemma~\ref{lemma U-acyclic}, we have $E^2_{i,j} = 0$ for all $j > 0$, so the
spectral sequence degenerates and
$H_i(G_\sigma, M) \cong E^2_{i,0} = H_i\big(\mathbb{Z}^{s-1}, M_{U_\sigma}\big)$, as desired.
\end{proof}

The flat dimension of $R$ as $R[G]$-module with trivial $G$-action is called the homological dimension of $G$ and denoted by $\mathrm{hd}_R(G)$ (cf. \cite[\S 4.4]{berihomologicaldimensiondisgroups}).
Equivalently, $\mathrm{hd}_R(G)$ is the smallest integer $n \geq 0$ such that, for each $R[G]$-module $M$, we have $H_r(G,M) = \lbrace 0 \rbrace$ for all $r>n$.
The number $\mathrm{hd}_R(G)$ depends on the commensurability class of $G$.
Indeed, for any two commensurable groups $H$ and $G$, $\mathrm{hd}_R(G)=\mathrm{hd}_R(H)$ whenever $\mathrm{hd}_R(G),\mathrm{hd}_R(H)< \infty$ (cf. \cite[Cor. 5.10]{berihomologicaldimensiondisgroups}).
The homological dimension can be $\infty$ in certain precise cases.
This occurs, for instance, when $G$ contains elements whose orders are not invertible in $R$.
When $\mathrm{hd}_R(G)=\infty$ and $G$ contains a finite-index subgroup $H$ such that $\mathrm{hd}_R(H)<\infty$, we can define the virtual $R$-homological dimension $\mathrm{vhd}_R(G)$ as $\mathrm{hd}_R(H)$.
This invariant is well-defined according to the discussion given above. The next result is then equivalent to Th. \ref{main teo HD}.

\begin{corollary}
One has $\mathrm{vhd}_R(\mathrm{SL}_2(\OS)) =s$, whenever $p \in R^{*}$.
\end{corollary}

\begin{proof}
Let $G$ be a finite-index subgroup of $\mathrm{SL}_2(\OS)$ whose torsion is $p$-primary. We can consider, for example, $G$ as a principal congruence subgroup of $\mathrm{SL}_2(\OS)$.
We will prove that $\mathrm{hd}_R(G) = s$.
Recall that $\mathrm{hd}_R(H) \leq \mathrm{hd}_R(G)$, for each subgroup $H \leq G$ (cf. \cite[Prop. 4.9]{berihomologicaldimensiondisgroups}).
Let $T_G \subset G$ be the subgroup isomorphic to $\mathbb{Z}^{s}$ given by Lemma \ref{lemma torus rank s}.
Then, we have $\mathrm{hd}_R(T_G)=\mathrm{hd}_R(\mathbb{Z}^{s})=s$, whence $\mathrm{hd}_R(G) \geq s$.
Now, it follows from Cor. \ref{cor direct sum} that, for any $r >s$ and any $R[G]$-module $M$, we have $H_{r}(G,M) \cong \bigoplus_{\sigma \in \Sigma_0} H_r(G_{\sigma},M)$.
Since for each $\sigma \in \Sigma_0$, we have $H_*(G_{\sigma},M) \cong H_*(\mathbb{Z}^{s-1},M_{U_{\sigma}})$ by Corollary \ref{coro hochschild serre} and $\mathrm{hd}_R(\mathbb{Z}^{\alpha})=\alpha$, we have $H_r(G_{\sigma},M) =\lbrace 0 \rbrace$, for all $r>s$. 
Hence, we have $H_{r}(G,M) =\lbrace 0 \rbrace$, for all $r>s$, whence $\mathrm{hd}_R(G) \leq s$.
\end{proof}

We now are able to prove Th. \ref{main theo intro pd of delta}.

\begin{proof}[Proof of Th. \ref{main theo intro pd of delta}]
If $s = 1$, then $U_\sigma = G_\sigma$ by Lemma~\ref{lemma structure Gs}, so $\Delta_{\mathfrak{U}} = \mathrm{St}$,
and the result is Theorem~\ref{theorem pd}. Assume that $s \geq 2$.
Applying $\mathrm{Tor}^{R[G]}_*(-, M)$ to
$0 \to \Delta_{\mathfrak{U}} \to \bigoplus_{\sigma} R[G/U_\sigma] \to R \to 0$
and using Shapiro's lemma, we obtain the long exact sequence
\begin{equation}\label{eq exact sequence of U_i}
\cdots \to \bigoplus_{\sigma \in \Sigma_0} H_i(U_\sigma, M) \to H_i(G, M)
\to H_i(G \bmod U_\sigma, M) \to \bigoplus_{\sigma \in \Sigma_0} H_{i-1}(U_\sigma, M)
\to \cdots .
\end{equation}
By Lemma~\ref{lemma U-acyclic}  the long
exact sequence~\eqref{eq exact sequence of U_i} gives, for every $n \geq 1$ and
every $R[G]$-module $M$,
\[
\mathrm{Tor}^{R[G]}_n(\Delta_{\mathfrak{U}}, M) = H_{n+1}(G \bmod U_\sigma, M) \cong H_{n+1}(G, M).
\]
Since $G$ has homological dimension $s$, the flat dimension of $\Delta_{\mathfrak{U}}$ is exactly $s-1$.
\end{proof}

The next result is immediate:

\begin{corollary}
Assume that $p\in R^{*}$ and that the torsion of $G$ is $p$-primary.
Then, for every $r>s$ and every $R[G]$-module $M$, we have $H_r(G \bmod U_{\sigma},M)=\lbrace 0\rbrace$. \qed
\end{corollary}

%\begin{remark}
%It is also true that for each $r >s$ and every $R[G]$-module $M$, we have $H_{r}(G,M) \cong \bigoplus_{\sigma \in \Sigma_0} H_r(U_{\sigma},M)$. 
%However, this equality is trivial since all the involved groups are zero.
%\end{remark}

\begin{comment}
Let $r > s$. By Corollary~5.2, we have
$H_r(G, M) \cong \bigoplus_{\sigma \in \Sigma_0} H_r(G_\sigma, M)$. By Corollary~6.3,
$H_r(G_\sigma, M) \cong H_r(\mathbb{Z}^{s-1}, M_{U_\sigma})$, which vanishes since the
Koszul resolution of $R$ over $R[\mathbb{Z}^{s-1}]$ has length $s-1$. Hence
$\mathrm{Tor}^{R[G]}_n(\Delta_{\mathfrak{U}}, M) = 0$ for all $n \geq s$.
This implies that $\mathrm{fd}_{R[G]}(\Delta_{\mathfrak{U}}) \leq s-1$.

Now, let $T' \subseteq G$ be a subgroup isomorphic to $\mathbb{Z}^{s}$,
given by Lemma~\ref{lemma torus rank s}. For each $g \in G$ and $\sigma \in \Sigma_0$,
the group $T' \cap gU_\sigma g^{-1}$ is a subgroup of an elementary abelian
$p$-group. Hence Mackey's formula and the argument of Lemma~6.2 give
$H_n\big(T', R[G/U_\sigma]\big) = 0$ for all $n \geq 1$. Applying $H_*(T', -)$ to
the short exact sequence
$$
0 \to \Delta_{\mathfrak{U}} \to \bigoplus_{\sigma \in \Sigma_0} R[G/U_\sigma] \to R \to 0,
$$
we obtain $H_{s-1}(T', \Delta_{\mathfrak{U}}) \cong H_s(T', R) \cong R$, by Koszul's
lemma. Therefore, by Shapiro's lemma,
\[
\mathrm{Tor}^{R[G]}_{s-1}\big(\Delta_{\mathfrak{U}}, R \otimes_{R[T']} R[G]\big)
\cong H_{s-1}(T', \Delta_{\mathfrak{U}}) \neq 0,
\]
and so $\mathrm{fd}_{R[G]}(\Delta_{\mathfrak{U}}) \geq s-1$.
\end{comment}

As in \S \ref{section Rel. Homology}, we now focus on describing the relative homology groups $H_r(G \bmod U_\sigma, M)$, for $r\leq s$ in terms of the Euler-Poincar\'e characteristic $\chi(G)$.
In particular, we will prove Th. \ref{main thm euler-unipotent for U}.
As in \S \ref{section Rel. Homology}, for each principal ideal domain
$D$ and each $D$-module $A$, we write $\mathrm{rk}_{D}(A):=\dim_{F}(A\otimes_{D}F)
\in \mathbb{Z}_{\geq0}\cup\lbrace\infty\rbrace$, where $F=\mathrm{Quot}(D)$.
Recall that $\mathrm{rk}_D$ is additive on short exact sequences. Hence the
alternating sum of the ranks of the terms of a bounded complex of $D$-modules of
finite rank equals that of its homology, and it vanishes for a bounded exact
sequence of $D$-modules of finite rank.

\begin{lemma}\label{lemma Koszul}
Let $D$ be a principal ideal domain, let $r \geq 1$, and let $N$ be a
$D[\mathbb{Z}^{r}]$-module of finite rank. Then $H_i(\mathbb{Z}^{r}, N)$
has finite rank for every $i$, vanishes for $i > r$, and
$$
\sum_{i=0}^{r} (-1)^i \mathrm{rk}_D H_i(\mathbb{Z}^{r}, N) = 0 .
$$
\end{lemma}

\begin{proof}
The Koszul resolution of the trivial module $D$ over $D[\mathbb{Z}^{r}]$
is a free resolution of length $r$, whose $i$-th term has rank $\binom{r}{i}$.
Tensoring it with $N$, we see that $H_i(\mathbb{Z}^{r}, N)$ is the $i$-th homology
group of a complex of the form
$$
0 \to N^{\binom{r}{r}} \to \cdots \to N^{\binom{r}{1}} \to N \to 0 .
$$
Hence $H_i(\mathbb{Z}^{r}, N)$ has finite rank and vanishes for $i > r$.
Moreover, by the additivity of $\mathrm{rk}_D$, we have
$$
\sum_{i=0}^{r} (-1)^i \mathrm{rk}_D H_i(\mathbb{Z}^{r}, N)
= \sum_{i=0}^{r} (-1)^i \binom{r}{i} \mathrm{rk}_D N
= (1-1)^{r} \mathrm{rk}_D N = 0.
$$

\end{proof}

\begin{corollary}\label{coro Euler-Gsigma}
Assume that $s \geq 2$ and that the torsion of $G$ is $p$-primary. Let
$D$ be a principal ideal domain with $p \in D^{*}$, and let $M$ be a
$D[G_\sigma]$-module of finite rank. Then $H_i(G_\sigma, M)$
has finite rank for every $i$, vanishes for $i \geq s$, and
$$
\sum_{i \geq 0} (-1)^i \mathrm{rk}_DH_i(G_\sigma, M) = 0 .$$
\end{corollary}

\begin{proof}
By Corollary~\ref{coro hochschild serre}, we have
$H_i(G_\sigma, M) \cong H_i\big(\mathbb{Z}^{s-1}, M_{U_\sigma}\big)$, and $M_{U_\sigma}$
has finite rank since it is a quotient of $M$. The result follows from
Lemma~\ref{lemma Koszul} with $r = s-1 \geq 1$.
\end{proof}

Now we prove Th. \ref{main thm euler-unipotent for U} which describes $H_r(G \bmod G_{\sigma},M)$, $r\leq s$, in terms of $\chi(G)$.

\begin{proof}[Proof of Th. \ref{main thm euler-unipotent for U}]
Since $p \in D^{*}$, we apply Lemma~\ref{lemma U-acyclic},
Corollary~\ref{coro hochschild serre} and Corollary~\ref{coro Euler-Gsigma} with
$R = D$.
Let $\Delta_{\mathfrak{H}} = \ker\big(\bigoplus_{\sigma \in \Sigma_0} D[G/H_\sigma]
\to D\big)$ be the kernel of the augmentation map with respect to a family of
subgroups $\mathfrak{H}=\lbrace H_{\sigma}\rbrace_{\sigma \in \Sigma_0}$ of $G$.
For instance, if $\mathfrak{H}=\mathfrak{G}$, then $\Delta_{\mathfrak{H}}=\mathrm{St}$.
Applying $\mathrm{Tor}^{D[G]}_*(-, M)$ to
$0 \to \Delta_{\mathfrak{H}} \to \bigoplus_{\sigma} D[G/H_\sigma] \to D \to 0$
and using Shapiro's lemma, we obtain the long exact sequence
\begin{equation}\label{eq exact sequence of H_i}
\cdots \to \bigoplus_{\sigma \in \Sigma_0} H_i(H_\sigma, M) \to H_i(G, M)
\to H_i(G \bmod H_\sigma, M) \to \bigoplus_{\sigma \in \Sigma_0} H_{i-1}(H_\sigma, M)
\to \cdots .
\end{equation}
By the additivity of $\mathrm{rk}_D$, if all the terms of
\eqref{eq exact sequence of H_i} have finite rank and vanish in large degrees, then
\begin{equation}\label{eq:Euler-LES}
\sum_{i \geq 1} (-1)^i \mathrm{rk}_D H_i(G \bmod H_\sigma, M)
= \sum_{i \geq 0} (-1)^i \mathrm{rk}_D H_i(G, M)
- \sum_{\sigma \in \Sigma_0} \sum_{i \geq 0} (-1)^i \mathrm{rk}_D H_i(H_\sigma, M).
\end{equation}

We apply \eqref{eq exact sequence of H_i} to the family
$\{G_\sigma\}_{\sigma \in \Sigma_0}$, for which $\Delta_{\mathfrak{H}} = \mathrm{St}$ by
Lemma~2.5. In this sequence, each $H_i(G, M)$ sits between
$\bigoplus_\sigma H_i(G_\sigma, M)$ and $H_i(G \bmod G_\sigma, M)$. Since $\Sigma_0$ is
finite, the former has finite rank and vanishes for $i \geq s$ by
Corollary~\ref{coro Euler-Gsigma}. The latter has finite rank and vanishes for
$i > s$ by Theorem~\ref{thm euler relative}~(a). Hence $H_i(G, M)$ has finite rank
and vanishes for $i > s$. Moreover, Equation~\eqref{eq:Euler-LES}, together with
Corollary~\ref{coro Euler-Gsigma} and Theorem~\ref{thm euler relative}~(b), gives
\begin{align*}
\sum_{i \geq 0} (-1)^i \mathrm{rk}_D H_i(G, M)
& = \sum_{i \geq 1} (-1)^i \mathrm{rk}_D H_i(G \bmod G_\sigma, M)
+ \sum_{\sigma \in \Sigma_0} \sum_{i \geq 0} (-1)^i \mathrm{rk}_D H_i(G_\sigma, M)\\
&= \chi(G) \cdot m + 0 .
\end{align*}

We now apply \eqref{eq exact sequence of H_i} to the family
$\mathfrak{U} = \{U_\sigma\}_{\sigma \in \Sigma_0}$.
By Lemma~\ref{lemma U-acyclic}, we have $H_j(U_\sigma, M) = 0$ for $j \geq 1$ and
$H_0(U_\sigma, M) = M_{U_\sigma}$. Thus the long exact sequence
\eqref{eq exact sequence of H_i} reduces to isomorphisms
$$
H_i(G \bmod U_\sigma, M) \cong H_i(G, M), \quad \forall i \geq 2,
$$
together with the exact sequence
$$
0 \to H_1(G, M) \to H_1(G \bmod U_\sigma, M)
\to \bigoplus_{\sigma \in \Sigma_0} M_{U_\sigma} \to M_G \to 0 .
$$
Since $M_{U_\sigma}$ and $M_G$ are quotients of $M$, they have finite rank.
Statement~(a) then follows. Finally, Equation~\eqref{eq:Euler-LES} gives
$$
\sum_{i \geq 1} (-1)^i \mathrm{rk}_D H_i(G \bmod U_\sigma, M)
= \sum_{i \geq 0} (-1)^i \mathrm{rk}_D H_i(G, M)
- \sum_{\sigma \in \Sigma_0} \mathrm{rk}_D M_{U_\sigma}
= \chi(G) \cdot m - \sum_{\sigma \in \Sigma_0} \mathrm{rk}_D M_{U_\sigma},
$$
which is Statement~(b).
\end{proof}

\begin{corollary}\label{coro euler-unipotent-trivial}
Assume that $s \geq 2$, and let $D$ be as in Theorem~\ref{main thm euler-unipotent for U}.
Taking $M = D$ with the trivial $G$-action, we have
\[
\sum_{i=1}^{s} (-1)^i \, \mathrm{rk}_D H_i(G \bmod U_\sigma, D) = \chi(G) - h,
\]
where $h=\# \Sigma_0$ is the cusp number $\# G\backslash \mathbb{P}^1(k)$.
\end{corollary}

\begin{remark}
For $s = 1$, Lemma~\ref{lemma structure Gs} gives $G_\sigma = U_\sigma$, and
Theorem~\ref{thm euler relative} already computes $H_*(G \bmod U_\sigma, M)$. The
correction term $\sum_\sigma \mathrm{rk}_D(M_{U_\sigma})$ in
Theorem~\ref{main thm euler-unipotent for U} appears for $s \geq 2$ because the torus
$\mathbb{Z}^{s-1}$ of $G_\sigma$ has vanishing Euler characteristic. In particular,
for $s \geq 2$ the identity in Conjecture~\ref{conj serre rel hom}~(c) holds for the
family $\{G_\sigma\}_{\sigma \in \Sigma_0}$, but not for
$\{U_\sigma\}_{\sigma \in \Sigma_0}$ unless $\sum_\sigma \mathrm{rk}_D(M_{U_\sigma}) = 0$.
\end{remark}

\section{Examples}\label{section exs}

Let us assume that $\OS$ is a principal ideal domain, so that $\mathrm{Pic}(\OS)=\lbrace e \rbrace$.
Then, the arithmetic group $\mathrm{SL}_2(\OS)$ acts transitively on $\mathbb{P}^1(k)$ by \cite[Th. 3.1]{BravoLoisel2}.
For each $\sigma \in \mathbb{P}^1(k)$, fix $\gamma=\gamma_{\sigma} \in \mathrm{SL}_2(\OS)$ so that $\sigma= \gamma \cdot \infty$.
Therefore, for each $G \subseteq \mathrm{SL}_2(\OS)$ we have
$$G_{\sigma}= \mathrm{Stab}_G(\sigma)= G \cap \mathrm{Stab}_{\mathrm{SL}_2(\OS)}(\sigma)=  G \cap \gamma \,\mathrm{Stab}_{\mathrm{SL}_2(\OS)}(\infty) \, \gamma^{-1}. $$
Now, assume that $G$ is a normal subgroup of $\mathrm{SL}_2(\OS)$. Then, we have
$$G_{\sigma}=  \gamma \, \big(G \cap \mathrm{Stab}_{\mathrm{SL}_2(\OS)}(\infty)\big) \, \gamma^{-1} =  \gamma \,  \mathrm{Stab}_{G}(\infty) \, \gamma^{-1}= \gamma G_{\infty} \gamma^{-1}.$$
Hence, the stabilizer $G_{\sigma}$ is isomorphic to $\mathrm{Stab}_{G}(\infty) = G \cap \mathrm{B}(k)$.
As in \S \ref{section Rel. Homology} we write $h=\# \Sigma_0$.
The next result is a direct consequence of Cor. \ref{cor direct sum}.

\begin{corollary}
Let us denote by $G^{\natural}$ the group of upper triangular matrices $G^{\natural}:=G \cap \mathrm{B}(k)$. For each $r >s$ and every $R[G]$-module $M$, we have $H_{r}(G,M) \cong H_r(G^{\natural},M)^{h}$. For $r=s$, we have the exact sequence
$0 \to H_s(G^{\natural},M)^{h} \to H_s (G,M) \to H_s(G \text{ mod } G_{\sigma},M).$ \qed 
\end{corollary}

Let $I$ be an ideal of $\OS$. The principal congruence subgroup $\Gamma_I$ of $\Gamma:=\mathrm{SL}_2(\OS)$ defined from $I$ is $\Gamma_I=\mathrm{ker}\left( \mathrm{SL}_2(\OS) \to  \mathrm{SL}_2(\OS/I) \right)$. It follows from \cite[Lemma 3.3]{MasonPremetSuryZalesskii} that the torsion of $\Gamma_I$ is $p$-primary.
Moreover, $\Gamma_I$ is a normal subgroup of $\mathrm{SL}_2(\OS)$.
It is not hard to see that $\Gamma_I^{\natural}=\Gamma_I \cap \mathrm{B}(k)$ is isomorphic to $ I \rtimes  \left((1+I) \cap \OS^{*}\right)$.
Since $\mathbb{F}$ is finite, we have $\OS/I$ is finite, whence $(1+I) \cap \OS^{*}$ is a finite-index subgroup of $\OS^{*}$, because it is the kernel of the group homomorphism $\OS^{*} \to (\OS/I)^{*}$.
It follows from Dirichlet's unit theorem that $(1+I) \cap \OS^{*}\cong \mathbb{Z}^{s-1}$.
Therefore, the next result follows.

\begin{corollary}
Let $\Gamma_I$ be a principal congruence subgroup of $\mathrm{SL}_2(\OS)$.
For each $r >s$ and every $R[\Gamma_I]$-module $M$, we have $H_{r}(\Gamma_I,M) \cong H_r( \mathbb{Z}^{s-1} \rtimes I,M)^{h}$. For $r=s$, we have the exact sequence
$0 \to H_s( \mathbb{Z}^{s-1} \rtimes I ,M)^{h} \to H_s (\Gamma_I,M) \to H_s(\Gamma_I \text{ mod } (\Gamma_I)_{\sigma},M).$ \qed 
\end{corollary}

\section*{Acknowledgments}
M.A. was supported by ANID Fondecyt postdoctoral grant 3261344.

C.B. was partially supported by ANID Fondecyt Initiation Grant No. 11260422.

\bibliographystyle{amsalpha}
\bibliography{refs.bib}

%\end{thebibliography}
\end{document}